\documentclass[11pt]{article}

\usepackage{amsmath,amsthm}
\usepackage {latexsym}
\usepackage{amssymb}

\newcommand{\beeq}{\begin{equation}}
\newcommand{\eneq}{\end{equation}}

\newcommand{\bear}{\begin{eqnarray}}
\newcommand{\eear}{\end{eqnarray}}

\newcommand{\supp}{\mbox{\rm supp}}

\newcommand{\half}{\frac{1}{2}}

\newcommand{\eps}{{\varepsilon}}
\newcommand{\R}{{\mathbb R}}

\newcommand{\Compl}{{\mathbb C}}

\newcommand{\bdd}{{\mathcal B}}

\newcommand{\calS}{{\mathcal S}}

\newcommand{\les}{\lesssim}

\newcommand{\Laplace}{\triangle}

\newcommand{\kato}{{\mathcal K}}

\newcommand{\la}{\lambda}

\def\nn{\nonumber}
\def\Ran{{\rm Ran}}

\def\Pac{P_{ac}}
\def\la{\langle}
\def\ra{\rangle}

\def\norm[#1][#2]{\|#1\|_{#2}}

\def\x{{\bf x}}

\def\norm[#1][#2]{\Vert #1 \Vert_{#2}}
\def\trace{{\rm trace}}
\def\tchi{\widetilde{\chi}}
\newcommand{\Hmm}[1]{\leavevmode{\marginpar{\tiny%
$\hbox to 0mm{\hspace*{-0.5mm}$\leftarrow$\hss}%
\vcenter{\vrule depth 0.1mm height 0.1mm width \the\marginparwidth}%
\hbox to 0mm{\hss$\rightarrow$\hspace*{-0.5mm}}$\\\relax\raggedright #1}}}

\newtheorem{theorem}{Theorem}
\newtheorem{lemma}[theorem]{Lemma}
\newtheorem{defi}[theorem]{Definition}
\newtheorem{cor}[theorem]{Corollary}
\newtheorem{prop}[theorem]{Proposition}

\theoremstyle{remark}

\begin{document}

\title{Dispersive estimates for Schr\"odinger operators in
dimension two}
\date{}

\author{W.\ Schlag\thanks{The author was partially supported by the NSF grant
DMS-0300081 and a Sloan Fellowship}}

\maketitle

\section{Introduction}

The purpose of this paper is to prove the following
result.

\begin{theorem}
\label{thm:main}
Let $V:\R^2\to\R$ be a measurable function
such that $|V(x)|\le C(1+|x|)^{-\beta}$, $\beta>3$. Assume in addition
that zero is a regular point of the spectrum of $H=-\Laplace+V$. Then
\[ \big\|e^{itH}P_{ac}(H) f\big \|_{\infty} \le C|t|^{-1}\|f\|_1\]
for all $f\in L^1(\R^2)$.
\end{theorem}

The definition of zero being a regular point amounts to the following, see
Jensen, Nenciu~\cite{JenNen} and Definition~\ref{defi:zero} below:
{\em Let $V\not\equiv 0$ and set $U={\mathrm sign}\, V$, $v=|V|^{\half}$.
Let $P_v$ be the orthogonal projection onto $v$ and set $Q=I-P_v$. Finally, let
\[ (G_0 f)(x):= -\frac{1}{2\pi}\int_{\R^2} \log|x-y|\, f(y)\, dy.\]
Then zero is regular iff $Q(U+vG_0v)Q$ is invertible on $QL^2(\R^2)$.}

Jensen and Nenciu study $\ker[Q(U+vG_0v)Q]$ on $QL^2(\R^2)$.
It can be completely described in terms of
solutions~$\Psi$ of $H\Psi=0$. In particular, its dimension is at most three plus the dimension of the
zero energy eigenspace, see Theorem~6.2 and Lemma~6.4  in~\cite{JenNen}.  The extra three
dimensions here are called resonances.
Hence, the requirement that zero is a regular point is the analogue of the usual
condition that zero is neither an eigenvalue nor a resonance of~$H$.

As far as the spectral properties of $H$ are concerned, we note that under the hypotheses
of Theorem~\ref{thm:main} the spectrum of $H$ on $[0,\infty)$ is purely absolutely
continuous, and that the spectrum is pure point on $(-\infty,0)$ with at most finitely
many eigenvalues of finite multiplicities. The latter follows for example from Stoiciu~\cite{mihai},
who obtained Birman-Schwinger type bounds in the case of two dimensions.

Theorem~\ref{thm:main} appears to be the first $L^1\to L^\infty$ bound
with $|t|^{-1}$ decay in~$\R^2$. Yajima~\cite{Yaj} and Jensen, Yajima~\cite{JenYaj}
proved the $L^p(\R^2)$ boundedness of the wave operators under stronger decay
assumptions on~$V(x)$, but only for $1<p<\infty$. Hence their result does
not imply Theorem~\ref{thm:main}. Local $L^2$ decay was studied by
Murata~\cite{Mur}, but he does not consider $L^1\to L^\infty$ estimates.

The first $L^1(\R^n)\to L^\infty(\R^n)$ bounds for $e^{itH}$ with $|t|^{-\frac{n}{2}}$ decay
were obtained by Journ\'e, Soffer, and Sogge~\cite{JSS}. However, their argument
 depends on the fact that $t^{-\frac{n}{2}}$ is integrable at $t=\infty$, and
thus only applies for $n\ge3$. In dimension $n=1$  Weder~\cite{Wed}
obtained the $|t|^{-\half}$-decay under some conditions on $V$ which were then relaxed by
Goldberg and the author~\cite{GolSch}. However, the case $n=2$ remained open.

As usual, the proof of Theorem~\ref{thm:main} breaks up into two regimes:
energies bigger than $\lambda_1$ and energies in $(0,\lambda_1)$.
Here $\lambda_1>0$ is some small constant.
The corresponding statements are Propositions~\ref{prop:high} and~\ref{prop:low_decay} below.
Theorem~\ref{thm:main} then follows by combining these two Propositions.
For energies in $(0,\lambda_1)$ we use the recent results of Jensen and Nenciu~\cite{JenNen}
on expansions of the resolvent $(H-(\lambda^2\pm i0))^{-1}$ for $\lambda$ close to zero.
Since we require somewhat finer estimates on various error terms, we give a complete
derivation of this expansion. However, we emphasize that this derivation is of course
merely a variant of a special case of the expansions in~\cite{JenNen}. In fact, the
main achievement of Jensen and Nenciu is to determine the expansion of the perturbed
resolvent in the presence of resonances and eigenvalues at zero.

\section{Energies separated from zero}
\label{sec:high}

The main purpose of this section is to prove the dispersive
estimate for the evolution restricted to energies
$[\lambda_1,\infty)$, $\lambda_1>0$. This will be accomplished by
an expansion of the perturbed resolvent into a finite Born series,
see~\eqref{eq:finsum} and~\eqref{eq:rest}. The main difficulty is
to obtain the dispersive bound for each term of the Born series.
This is done in Lemma~\ref{lem:split} below. For the
remainder~\eqref{eq:rest} in the Born expansion, which still
contains the perturbed resolvent, we use the limiting absorption
principle. The approach in this section is modelled after that
in~\cite{GolSch}, which in turn had its origins in the work of
Rodnianski and the author~\cite{RodSch}.

\noindent Lemma~\ref{lem:sp} is a variant of the standard stationary phase method.
In what follows, the notation $x\les y $ means that $x\le Cy$ for some constant $C$.

\begin{lemma}
\label{lem:sp}
Let $\phi(0)=\phi'(0)=0$ and $1\le \phi''\le C$. Then
\beeq
\label{eq:delta}
\left| \int_{-\infty}^\infty e^{it\phi(x)}\, a(x)\,dx\right|
\les \delta^2 \int \Big(\frac{|a(x)|}{\delta^2+|x|^2} + \chi_{[|x|>\delta]}\;\frac{|a'(x)|}{|x|}\Big)\,dx
\eneq
where $\delta=|t|^{-\half}$.
\end{lemma}
\begin{proof} With $\eta$ being a standard cut-off one has
\bear
\Bigl| \int_{-\infty}^\infty e^{it\phi(x)} a(x)\,dx \Bigr| &\le&
\Bigl| \int_{-\infty}^\infty e^{it\phi(x)} a(x)\eta(x/\delta)\,dx \Bigr| \nn \\
&& \qquad\qquad + \Bigl| \int_{-\infty}^\infty e^{it\phi(x)} a(x)(1-\eta(x/\delta))\,dx \Bigr| \nn \\
&\les& \int_{|x|<\delta} |a(x)|\,dx +
\delta^2 \int \Bigl| \Bigl(\frac{a(x)(1-\eta(x/\delta))}{\phi'(x)}\Bigr)'\Bigr|\, dx \nn \\
&\les& \int_{|x|<\delta} |a(x)|\,dx +
\delta^2 \int_{|x|\gtrsim \delta} \Bigl| \Bigl(\frac{|a(x)|}{|x|^2}+\frac{|a'(x)|}{|x|}\Bigr)\Bigr|\, dx, \nn
\eear
as claimed.
\end{proof}

It is well-known that
\[  R_0^{\pm}(\lambda^2)(x,y)= (-\Delta-(\lambda^2\pm  i0))^{-1}(x,y)=\pm \frac{i}{4} H_0^{\pm}(\lambda|x-y|), \]
where $H_0^{\pm}$ are the Hankel functions of order zero with $H_0^{-}=\overline{H_0^{+}}$.
 They have the form
\[ H_0^{+}(y)=e^{i(y-1)}\omega(y)\chi_{[y>1]}+\omega(y)\chi_{[0<y<1]}\]
 and satisfy the bounds
 $|\omega(y)|\les |y|^{-\half}$ if
$y\gtrsim 1$ and $|\omega(y)|\les |\log y|$ of $0<y<\half$.
Moreover, one has for all positive integers~$\nu$,
\bear
 |\omega^{(\nu)}(y)| &\les& |y|^{-\nu-\half} \text{\ \ if \ \ }y\gtrsim 1 \nn \\
 |\omega^{(\nu)}(y)| &\les& |y|^{-\nu} \text{\ \ if \ \ }0<y< 1. \nn
\eear
Set $\omega_+(y)=\chi_1(y/y_0)\omega(y)$ where
$\chi_1(y)=0$ if $y\le 1$ and $=1$ if $y\ge 2$. Here $y_0\gg1$ is a fixed constant.
Define $\omega_{-}(y)$ via $\omega=\omega_{+}+\omega_{-}$, i.e.,
$\omega_{-}(y)=(1-\chi_1(y/y_0))\omega(y)$ (in Section~\ref{sec:low} the functions
 $\omega_{+}$ and $\omega_{-}$
will take on a different meaning, not to be confused with the one here). Let
\beeq
\label{eq:kato} \|V\|_{\kato} :=  \sup_{x\in\R^2} \int_{\R^2} \Big(1+\log^{-}|x-y|\Big)^2\,|V(y)|\, dy,
\eneq
where $\log^{-}u=-\chi_{[0<u<1]}\,\log u$.
Finally, pick a cut-off $\chi_2$ so that
 $\chi_2(y)=1$ if $y\le1$ and $\chi_2(y)=0$ if $y\ge2$.
The following lemma is one of the two main technical ingredients
of the high energy part. We urge the reader not to be distracted
by the technical appearance of the proof. Indeed, the
bound~\eqref{eq:osc} can be derived heuristically as an immediate
consequence of stationary phase. However, some cases do need to be
distinguished due to various cut-offs in the integrand.

\begin{lemma}
\label{lem:split}
Assume $\|V\|_{\kato} <\infty$. Let $\{1,2,\ldots,m\}=J\cup J^*$ be a partition. Then
\begin{align}
& \sup_{\substack{L\ge1 \\x_0,x_{m}\in \R^2}} \int_{\R^{2(m-1)}}
\Bigl| \int_0^\infty \lambda\,e^{i(t\lambda^2\pm \lambda\sum_{j\in J} |x_{j-1}-x_j|)}
\chi_1(\lambda)\chi_2(\lambda/L) \prod_{j\in J} \omega_{+}(\lambda|x_{j-1}-x_j|)\nn \\
& \qquad  \prod_{\ell\in J^*}
\omega_{-}(\lambda|x_{\ell-1}-x_\ell|) \, d\lambda \Bigr|\;\prod_{k=1}^{m-1} |V(x_k)| \, dx_1\ldots dx_{m-1} \les |t|^{-1} \|V\|_{\kato}^{m-1}\label{eq:osc}
\end{align}
with a constant that only depends on $m$.
\end{lemma}
\begin{proof} The heuristic reason for this bound is as follows:
Let $d_j=|x_{j}-x_{j-1}|$ and $s=\sum_{j\in J}d_j$. If there is a critical point of the phase,
then it is $\lambda_0=\frac{s}{2t}$ (assuming $t>0$). We may assume that $\lambda_0 \gtrsim 1$, otherwise
the integrand vanishes at $\lambda_0$. Using stationary phase, the inner integral is then bounded by
\[ t^{-\half} \lambda_0 \, (\lambda_0 s)^{-\half} \, \prod_{\ell\in J^*} \log^{-}(\lambda_0 d_\ell)
   \les  t^{-1}  \prod_{\ell\in J^*} \log^{-}(d_\ell).
\]
Inserting this bound into~\eqref{eq:osc} then yields the desired result by an application of
Cauchy-Schwartz, see~\eqref{eq:Acont} below.

To make this rigorous, we start off integrating by parts: Then
\bear
&& |t|\,\Bigl|\int_0^\infty \lambda\, e^{i(t\lambda^2\pm \lambda s)}  \chi_1(\lambda)\chi_2(\lambda/L)
\prod_{j\in J} \omega_{+}(\lambda d_j)\, \prod_{\ell\in J^*} \omega_{-}(\lambda d_\ell) \, d\lambda
\Bigr| \nn \\
&& \les \Bigl| \int_0^\infty e^{i(t\lambda^2\pm \lambda s)} \chi_1'(\lambda)\chi_2(\lambda/L)
\prod_{j\in J} \omega_{+}(\lambda d_j)\, \prod_{\ell\in J^*} \omega_{-}(\lambda d_\ell) \, d\lambda
\Bigr| \nn \\
&& + \frac{1}{L}\Bigl| \int_0^\infty e^{i(t\lambda^2\pm \lambda s)} \chi_1(\lambda)\chi_2'(\lambda/L)
\prod_{j\in J} \omega_{+}(\lambda d_j)\, \prod_{\ell\in J^*} \omega_{-}(\lambda d_\ell) \, d\lambda
\Bigr| \nn \\
&& + s\Bigl| \int_0^\infty e^{i(t\lambda^2\pm \lambda s)} \chi_1(\lambda)\chi_2(\lambda/L)
\prod_{j\in J} \omega_{+}(\lambda d_j)\, \prod_{\ell\in J^*} \omega_{-}(\lambda d_\ell) \, d\lambda
\Bigr| \label{eq:C-} \\
&& + \sum_{k\in J}d_k\Bigl| \int_0^\infty e^{i(t\lambda^2\pm \lambda s)} \chi_1(\lambda)\chi_2(\lambda/L)
\,\omega_{+}'(\lambda d_k)\prod_{\substack{j\in J\\j\ne k}} \omega_{+}(\lambda d_j)\, \prod_{\ell\in J^*} \omega_{-}(\lambda d_\ell) \, d\lambda
\Bigr| \nn \\
&& + \sum_{k\in J^*}d_k\Bigl| \int_0^\infty e^{i(t\lambda^2\pm \lambda s)} \chi_1(\lambda)\chi_2(\lambda/L)
\prod_{j\in J} \omega_{+}(\lambda d_j)\,\omega_{-}'(\lambda d_k) \prod_{\substack{\ell\in J^*\\ \ell\ne k}} \omega_{-}(\lambda d_\ell) \, d\lambda
\Bigr| \nn \\
&& =: A^{\pm}+B^{\pm}+C^{\pm}+D^\pm+E^\pm. \nn
\eear
Let $k(x,y):=1+\log^{-}|x-y|$. Then since $\log^{-}$ is decreasing,
\bear
|A^{\pm}| &\les& \int |\chi_1'(\lambda)|\, \prod_{j\in J^*} |\omega_{-}(\lambda d_j)|\, d\lambda \nn \\
&\les& \int |\chi_1'(\lambda)|\, \prod_{j\in J^*}(1+ \log^{-}(\lambda d_j))\, d\lambda \les
\int |\chi_1'(\lambda)|\, \prod_{j\in J^*}(1+ \log^{-}( d_j))\, d\lambda \nn \\
&\les& \prod_{j\in J^*} k(x_{j-1},x_j). \nn
\eear
Hence the contribution of $A^{\pm}$ to~\eqref{eq:osc} is
\bear
&\les& \int_{\R^{2(m-1)}} \prod_{j=1}^m k(x_{j-1},x_j) \prod_{j=1}^{m-1}
|V(x_j)|\, dx_1\ldots dx_{m-1}\nn \\
&\les& \int \Big(k^2(x_0,x_1)|V(x_1)|+|V(x_1)|k^2(x_1,x_2)\Big)\prod_{j=2}^{m-1} |V(x_j)|k(x_j,x_{j+1})
\, dx_1\ldots dx_{m-1} \nn \\
&\les& \|V\|_{\kato}^{m-1}. \label{eq:Acont}
\eear
For the remainder of the proof we set
\[P_* = \prod_{j\in J^*} k(x_{j-1},x_j)\]
with the understanding that $P_*=1$ if $J^*=\emptyset$.
Similarly, for $L\ge1$, one has that
\[
|B^{\pm}| \les \frac{1}{L}\int_0^\infty |\chi_2'(\lambda/L)| \prod_{j\in J^*} |\omega_{-}(\lambda d_j)|
\, d\lambda
\les P_*.
\]
Hence the contribution by $B^\pm$ to \eqref{eq:osc} is again $\les \|V\|_{\kato}^{m-1}$.
The terms $D^{\pm}, E^{\pm}$ are also easy to deal with. Indeed, one has
\begin{align*}
|D^{\pm}| &\les \sum_{k\in J}d_k \int_0^\infty
(1+\lambda d_k)^{-\frac32}\prod_{\substack{j\in J\\j\ne k}} (1+\lambda d_j)^{-\half}\, d\lambda\, P_*\\
&= \int_0^\infty \Big[-2\prod_{j\in J} (1+\lambda d_j)^{-\half}\Big]'\, d\lambda\;P_* = 2\,P_*.
\end{align*}
As far as $E^{\pm}$ is concerned, we conclude similarly that (with some small constant $c>0$)
\begin{align*}
|E^{\pm}| &\les \sum_{k\in J^*}d_k \int_1^\infty
(d_k \lambda)^{-1}\chi_{[c\lambda d_k\le 1]}\prod_{\substack{j\in J^*\\j\ne k}} \log^{-}(c\lambda d_j)\,
d\lambda \\
&= \int_1^\infty \Big[-\prod_{j\in J^*} \log^{-}(c\lambda d_j)\Big]'\, d\lambda \les P_*.
\end{align*}
We now apply Lemma~\ref{lem:sp} to
$C^-$ with $\phi_{-}(\lambda)=\lambda^2-\lambda\frac{s}{t}$ and
\beeq
\label{eq:adef}
a(\lambda)=a_{+}(\lambda)\,\prod_{j\in J^*} \omega_{-}(\lambda d_j), \quad
a_{+}(\lambda) = \chi_1(\lambda)\chi_2(\lambda/L) \,\prod_{j\in J} \omega_{+}(\lambda d_j).
\eneq
Note that
\bear
|a(\lambda)| &\les& a_{+}(\lambda) \prod_{j\in J^*} k(x_{j-1},x_j)=a_+(\lambda)P_* \label{eq:aab} \\
|a'(\lambda)| &\les& |a_{+}'(\lambda)| P_*
+ \sum_{j\in J^*} |a_{+}(\lambda)| \lambda^{-1}\chi_{[\lambda d_j \les 1]}
 \prod_{\substack{k\in J^*\\ k\ne j}} (1+ \log^{-}(\lambda d_k)).\label{eq:a'ab}
\eear
Set $g(y)=1+\log^{-}(y)$ so that $g'(y)=-\chi_{[0<y<1]}y^{-1}$. Define
\beeq
\label{eq:bdef}
 b(\lambda) = \prod_{j\in J^*} g(c\lambda d_j)\text{\ \ with some small\ \ }c>0.
\eneq
Then $0<b(\lambda) \les \prod_{j\in J^*} k(x_{j-1},x_j)$ for $\lambda\gtrsim 1$ and
\beeq
|a'(\lambda)| \les |a_{+}'(\lambda)| P_* 
+ |a_+(\lambda)| |b'(\lambda)|
\les |a_{+}'(\lambda)| P_* 
+ a_0(\lambda) (-b'(\lambda))\label{eq:a0b}
\eneq
where $a_0(\lambda) = \chi_1(\lambda)\chi_2(\lambda/L) \prod_{j\in J} (1+\lambda d_j)^{-\half}$.

\smallskip
\noindent One has $\phi_{-}'(\lambda_0)=0$ for $\lambda_0=\frac{s}{2t}$. We first assume
that $\chi_1(\lambda_0)\ne0$ as well as $\lambda_0\in \supp(\omega_{+}(d_j \cdot))$
for each $j\in J$. These assumptions translate into $\lambda_0\gtrsim 1$ and
$\min_{j\in J} \lambda_0 d_j\gg 1$. The latter condition implies that
$\lambda_0^2=\frac{s\lambda_0}{2t}\gg t^{-1}$ and thus $\lambda_0\gg \delta=|t|^{-\half}$.
Next, we use Lemma~\ref{lem:sp} to bound $|C^{-}|$. On the one hand,
\bear
&& s\delta^2 \int \frac{|a(\lambda)|P_*^{-1}}{\delta^2+|\lambda-\lambda_0|^2}\, d\lambda \nn \\
&\les&
s\delta^2 \int_{\lambda_0-\delta}^\infty \frac{(1+\lambda s)^{-\half}}{\delta^2+(\lambda-\lambda_0)^2}\, d\lambda
+ s\delta^2 \int_1^{\lambda_0-\delta} \frac{(1+\lambda s)^{-\half}}{(\lambda-\lambda_0)^2}\, d\lambda
\nn \\
&\les& s\delta^2 (1+\lambda_0 s)^{-\half} \int_{\lambda_0-\delta}^\infty
\frac{1}{\delta^2+(\lambda-\lambda_0)^{2}}\, d\lambda + \sqrt{s}\delta^2 \int_1^{\lambda_0-\delta} \frac{d\lambda}{\sqrt{\lambda}(\lambda_0-\lambda)^2} \nn \\&\les& \sqrt{s} \lambda_0^{-\half} \delta + \sqrt{s}\delta^2
\lambda_0^{-\frac32} \nn \\
&\les& \sqrt{\frac{s}{t}}\; \lambda_0^{-\half} + \sqrt{\frac{s}{t}}\; \lambda_0^{-\half}\delta \lambda_0^{-1}
\les 1. \label{eq:c-1}
\eear
On the other hand, see~\eqref{eq:a0b}, an integration by parts yields
\begin{align}
& s\delta^2 \int_{|\lambda-\lambda_0|>\delta} \frac{|a'(\lambda)|}{|\lambda-\lambda_0|}\, d\lambda  \nn \\
& \les s\delta^2\int_{|\lambda-\lambda_0|>\delta} \frac{|a_{+}'(\lambda)|P_*}{|\lambda-\lambda_0|}\,
d\lambda
+ s\delta^2\int_{|\lambda-\lambda_0|>\delta} \frac{a_0(\lambda)(-b'(\lambda))}{|\lambda-\lambda_0|}\,
d\lambda \nn \\
&\les s\delta^2\int_{|\lambda-\lambda_0|>\delta} \frac{|a_{+}'(\lambda)|P_*}{|\lambda-\lambda_0|}\,
d\lambda
+ s\delta^2 \int_{|\lambda-\lambda_0|>\delta} \frac{a_0'(\lambda) b(\lambda)}{|\lambda-\lambda_0|}\,
d\lambda \label{eq:12} \\
& +  s\delta^2 \int_1^{\lambda_0-\delta} \frac{a_0(\lambda) b(\lambda)}{(\lambda-\lambda_0)^2}\,
d\lambda  + s\delta^2  \frac{a_0(\lambda) b(\lambda)}{|\lambda-\lambda_0|}\Bigg|_{\lambda-\lambda_0=\pm \delta}. \label{eq:34}
\end{align}
By the estimates leading up to \eqref{eq:c-1} one has
$\eqref{eq:34}\les P_*$.  On the other hand,
\begin{align}
&\eqref{eq:12}
 \les s\delta^2P_* \int_{\lambda_0+\delta}^\infty \Bigl(\frac{1}{L}|\chi_2'(\lambda/L)|(1+\lambda s)^{-\half}
+\lambda^{-1}(1+\lambda s)^{-\half} \Bigr) \frac{d\lambda}{\lambda-\lambda_0} \nn \\
& + s\delta^2 P_* \int_1^{\lambda_0-\delta} \Bigl[ |\chi_1'(\lambda)|(1+s)^{-\half} +
\frac{|\chi_2'(\lambda/L)|}{L(1+\lambda s)^{\half}} \nn \\
&\qquad\qquad + \sum_{j\in J} d_j (1+\lambda d_j)^{-\frac32} \prod_{\substack{k\in J\\k\ne j}} (1+\lambda d_k)^{-\half} \Bigr] \frac{d\lambda}{\lambda_0-\lambda}. \label{eq:resum}
\end{align}
It will be convenient to resum the expression on the right-hand side of~\eqref{eq:resum} by rewriting
it as a derivative. This yields
\begin{align}
\eqref{eq:resum}
& \les \sqrt{s}\delta \lambda_0^{-\half} P_* + s\delta^2 P_* \int_1^{\lambda_0-\delta} \Big[-\prod_{k\in J}
(1+\lambda d_k)^{-\half} \Big]' \frac{d\lambda}{\lambda_0-\lambda} \nn \\
& \les P_*  + s\delta^2 P_* \int_1^{\lambda_0-\delta} \prod_{k\in J}
(1+\lambda d_k)^{-\half}  \frac{d\lambda}{(\lambda_0-\lambda)^2}  \nn \\
& + s\delta^2 P_*  \prod_{k\in J} (1+\lambda d_k)^{-\half}  \frac{1}{\lambda_0-\lambda} \Bigg|_{\lambda=1}\nn \\ & \les P_* + \sqrt{s}\delta^2 P_* \int_1^{\lambda_0-\delta}
  \frac{d\lambda}{\lambda^{\half}(\lambda_0-\lambda)^2} + s\delta^2 \lambda_0^{-1}P_* \nn \\
& \les P_* + (\sqrt{s} \delta^2\lambda_0^{-\frac32} + \sqrt{s}\delta \lambda_0^{-\half})P_* \les P_*. \nn
\end{align}
In view of the preceding, $|C^{-}|\les P_*$ provided $\lambda_0\gtrsim 1$ and
$\min_{j\in J}\lambda_0 d_j \gg 1$. This gives the desired contribution to~\eqref{eq:osc}.

\noindent
Now suppose that $\lambda_0\gtrsim 1$ but $\min_{j\in J}\lambda_0 d_j \les 1$. Let
$\mu=\min_{j\in J} d_j$ so that $\mu\les \lambda_0^{-1}$.
By construction, $\supp(a)\subset [C\mu^{-1},\infty)$ for some large $C$.
Therefore, $\lambda\gg\lambda_0$ as well as $\lambda-\lambda_0\gtrsim \lambda$ on $\supp(a)$.
By Lemma~\ref{lem:sp},
\begin{align*}
|C^{-}| &\les s\delta \max_{[\lambda_0-\delta,\lambda_0+\delta]} |a(\lambda)| + s\delta^2 \int_{\lambda_0+\delta}^\infty \Bigl[ \frac{|a(\lambda)|}{(\lambda-\lambda_0)^2} + \frac{|a'(\lambda)|}{\lambda-\lambda_0}\Bigr]\,d\lambda \\
& \les \sqrt{s}\,\delta \lambda_0^{-\half} \chi_{[\delta\ll \lambda_0]}P_* + s\delta
\big(1+s\mu^{-1}\big)^{-\half} \chi_{[\delta \gtrsim \lambda_0]}P_* + s\delta^2 \int_{\mu^{-1}}^\infty \Bigl[ \frac{|a(\lambda)|}{\lambda^2}+\frac{|a'(\lambda)|}{\lambda}\Bigr]\, d\lambda \\
&\les P_* + s\delta^2 \int_{\mu^{-1}}^\infty \Bigl[ \frac{|a(\lambda)|}{\lambda^2}+\frac{|a'(\lambda)|}{\lambda}\Bigr]\, d\lambda.
\end{align*}
To bound the integral we use
\begin{align*}
|a(\lambda)| &\les (1+s\lambda)^{-\half}\chi_{[\lambda>\mu^{-1}]} P_* \\
|a'(\lambda)| &\les \lambda^{-1}(1+s\lambda)^{-\half}\chi_{[\lambda>\mu^{-1}]} P_*,
\end{align*}
see \eqref{eq:aab} and~\eqref{eq:a'ab}. Therefore,
\begin{align}
s\delta^2 \int_{\mu^{-1}} \Bigl[ \frac{|a(\lambda)|}{\lambda^2}+\frac{|a'(\lambda)|}{\lambda}\Bigr]\, d\lambda &\les s\delta^2 P_* \int_{\mu^{-1}}^\infty \frac{(1+\lambda s)^{-\half}}{\lambda^2}\, d\lambda \label{eq:intmu}\\
&\les \sqrt{s}\delta^2 \mu^{\frac32} P_* \les s\delta^2 \mu P_* = \lambda_0\mu P_* \les P_*,\nn
\end{align}
where we used $\mu\leq s$ to pass to the second inequality in the second line.
It remains to consider the case when $\lambda_0\ll 1$.
Note that $a(\lambda)=0$ if $\min_{j\in J} \lambda d_j\les 1$, which is the same as
$\lambda \les \mu^{-1}$. Also, $a(\lambda)=0$ is $\lambda \le 1$.
Then, again via Lemma~\ref{lem:sp}, one obtains as in~\eqref{eq:intmu},
\begin{align*}
|C^{-}| &\les \sqrt{s}\delta^2 P_* \int_{1+\mu^{-1}}^\infty \lambda^{-\frac52}\, d\lambda \les \sqrt{s}\delta^2 (1+\mu^{-1})^{-\frac32}P_*\\
&\les \sqrt{s}\delta^2 \chi_{[\mu>1]}P_* + \sqrt{s}\delta^2\mu^{\frac32} \chi_{[\mu<1]}P_* \\
&\les \frac{s}{t}P_*+\frac{s}{t}\mu \chi_{[\mu<1]}P_* \les (\lambda_0 + \lambda_0\mu\chi_{[\mu<1]})P_* \les P_*.
\end{align*}
The lemma is proved.
\end{proof}

\begin{prop}\label{prop:high}
Assume that $|V(x)|\les (1+|x|)^{-\beta}$ for some  $\beta > 2$.
Let $H=-\Laplace+V$ and $\lambda_1>0$ be fixed. Then
\[ \sup_{L\ge 1}
\Big|\Big \la e^{itH} \chi_2(\sqrt{H}/L)\chi_1(\sqrt{H}/\lambda_1) \, f,g \Big\ra \Big|
\les |t|^{-1}\|f\|_1\|g\|_1
\]
for all $f,g\in \calS(\R^2)$. The constant here depends only on $V$ and $\lambda_1$.
\end{prop}
\begin{proof}
Let $R_V^{\pm}(\lambda^2)=(-\Laplace+V-(\lambda^2\pm i0))^{-1}$ be the perturbed resolvent.
It satisfies the limiting absorption principle, see Agmon~\cite{agmon},
\beeq
\label{eq:lim_ap}
 \|R_V^{\pm}(\lambda^2)\|_{L^{2,\sigma}(\R^2)\to L^{2,-\sigma}(\R^2)} < \infty
\eneq provided $\sigma>\half$. Here $L^{2,\sigma}(\R^2)$ is the
usual weighted space with norm
\[ \|f\|_{2,\sigma}=\Bigl(\int_{\R^2} (1+|x|)^{2\sigma}|f(x)|^2\, dx\Bigr)^{\half}.\]
In addition, one has
\[ \|\partial_\lambda R_V^{\pm}(\lambda^2)\|_{L^{2,\sigma}(\R^2)\to L^{2,-\sigma}(\R^2)} < \infty\]
provided $\sigma>\frac32$. The free resolvent satisfies the same bounds with some
decay in $\lambda$, say $\lambda^{-\alpha}$. The exact value of $\alpha>0$ is not relevant for our
purposes.
One has
\begin{align}
& \Big \la e^{itH} \chi_2(\sqrt{H}/L)\chi_1(\sqrt{H}/\lambda_1) \, f,g \Big\ra
\nn \\
& =  \int_0^\infty
e^{it\lambda^2}\lambda\, \chi_2(\lambda/L)\chi_1(\lambda/\lambda_1)\, \Big \la [R_V^{+}(\lambda^2)-R_V^{-}(\lambda^2)]f,g \Big\ra
\, \frac{d\lambda }{\pi i}. \label{eq:spec_theo_high}
\end{align}
We use the resolvent expansion
\begin{align}
 R_V^{\pm}(\lambda^2) & =  \sum_{\ell=0}^{2m+2} R_0^{\pm}(\lambda^2)(-VR_0^{\pm}(\lambda^2))^\ell \label{eq:finsum} \\
& + R_0^{\pm}(\lambda^2)(VR_0^{\pm}(\lambda^2))^mV R_V^{\pm}(\lambda^2)V(R_0^{\pm}(\lambda^2)V)^m
R_0^{\pm}(\lambda^2). \label{eq:rest}
\end{align}
Here $m$ is a positive integer that depends on~$\alpha$.
Recall that
\[ R_0^{\pm}(\lambda^2)(x,y) = \pm\frac{i}{4}\,H_0^{\pm}(\lambda |x-y|)\] (the Hankel functions
of order zero).
By Lemma~\ref{lem:split} each of
the finitely many terms in~\eqref{eq:finsum} leads to the desired time-decay in~\eqref{eq:spec_theo_high}.
In fact, this only requires that $\|V\|_{\kato}<\infty$.
For the term~\eqref{eq:rest} one proceeds as in the three-dimensional argument
via the limiting absorption principle and stationary phase, see~\cite{GolSch}.
Following Yajima~\cite{Yaj}, set
\[ G_{\pm,x}(\lambda)(x_1):= e^{\mp i\lambda|x|}R^{\pm}_0(\lambda^2)(x_1,x).\]
Removing $f,g$ from~\eqref{eq:spec_theo_high}, we are led to proving that
\begin{align}
& \Big|
\int_0^\infty e^{it\lambda^2}e^{\pm i\lambda(|x|+|y|)}\;\chi_2(\lambda/L)\,\chi_1(\lambda/\lambda_1)
\lambda \Big\la VR^{\pm}_V(\lambda^2)V (R_0^{\pm}(\lambda^2)V)^m G_{\pm,y}(\lambda), \nn \\
&\quad (R_0^{\mp}(\lambda^2)V)^m G_{\pm,x}^*(\lambda) \Big\ra \, d\lambda
\Big| \les |t|^{-1} \label{eq:main}
\end{align}
uniformly in $x,y\in\R^2$ and $L\ge 1$.
Next, we check that the derivatives of $G_{+,x}(\lambda)$ satisfy the
estimates (for $\lambda>\lambda_1>0$)
\begin{align}
\sup_{x\in\R^3} \Big\| \partial_\lambda^j G_{\pm,x}(\lambda)\Big\|_{L^{2,-\sigma}} &< C_{j,\sigma}\,\lambda^{-\half}\la x\ra^{-\eps} \text{\ \ provided\ \ } \sigma > \frac12 + j \label{eq:Gest1}\\
\sup_{x\in\R^3} \Big\| \partial_\lambda^j G_{\pm,x}(\lambda)\Big\|_{L^{2,-\sigma}} &< C_{j,\sigma}(\lambda {\la x\ra})^{-\half}
\text{\ \ provided\ \ } \sigma>1+j \label{eq:Gest2}
\end{align}
for all $j\ge 0$. The small $\eps>0$ in~\eqref{eq:Gest1} depends on~$\sigma$.
The bound~\eqref{eq:Gest2} is Lemma~3.1 in~\cite{Yaj}. Alternatively, both bounds
follow easily by writing
$H_0^{\pm}(u)=e^{\pm iu} \rho_{\pm}(u)$ where $|\rho_{\pm}(u)|\les |\log^{-}(u)|\chi_{[0<u<\half]}
+u^{-\half}\chi_{[u>\half]}$. Thus, consider
\begin{align*}
&\Big\|\partial_\lambda^j e^{\pm i\lambda(|y-x|-|x|)}\rho_{\pm}(\lambda|x-y|)\la y\ra^{-\sigma} \Big\|_{L^2_y(\R^2)}^2 \\
& \les \int \la y\ra^{2(j-\sigma)}| \rho_{\pm}(\lambda|x-y|)|^2\, dy \\
&  \les \int_{[\lambda|x-y|<\half]} \la y\ra^{2(j-\sigma)} |\log(\lambda|x-y|)|^2\, dy
+ \lambda^{-1}\int_{\R^2} \la y\ra^{2(j-\sigma)} |y-x|^{-1}\, dy \\
&\les \la x\ra^{2(j-\sigma)}\lambda^{-2} + \lambda^{-1} \la x\ra^{-1}\chi_{[\sigma>j+1]}
+\lambda^{-1} \la x\ra^{2(j-\sigma)+1}\chi_{[\sigma<j+1]}.
\end{align*}
The stated bounds now follow by making the appropriate choices of $\sigma$ depending on~$j$.

\noindent Rewrite the integral in~\eqref{eq:main} in the form (with $L=\infty$)
\begin{equation}
\label{eq:Ipm}
I^{\pm}(t,x,y):=\int_0^\infty e^{it\lambda^2 \pm i\lambda(|x|+|y|)} a^{\pm}_{x,y}(\lambda)\, d\lambda.
\end{equation}
By the aforementioned bounds on $R_0^{\pm}(\lambda^2)$ and $R_0^{\pm}(\lambda^2)$ on
weighted $L^2$-spaces, which provide decay in~$\lambda$, as well as~\eqref{eq:Gest1}, \eqref{eq:Gest2},
one concludes that $a^{\pm}_{x,y}(\lambda)$ has one derivative in~$\lambda$ and
\begin{align}
\Big|  a^{\pm}_{x,y}(\lambda)\Big| &\les (1+\lambda)^{-2} (\la x\ra\la y\ra)^{-\half}
 \text{\ \ for all\ \ }\lambda>\lambda_1  \label{eq:adec} \\
\Big|\partial_\lambda a^{\pm}_{x,y}(\lambda)\Big| &\les (1+\lambda)^{-2}\la x\ra^{-\eps} \qquad \text{\ \ for\ all\ \ }\lambda>\lambda_1, \label{eq:one_der}
\end{align}
which in particular justifies taking $L=\infty$ in~\eqref{eq:Ipm}.
This requires that one takes $m$ sufficiently large and
that $|V(x)|\les (1+|x|)^{-\beta}$ for some $\beta> 2$. The latter condition
arises as follows: Consider, \eqref{eq:adec}. Then by~\eqref{eq:Gest1} and the limiting
absorption principle, respectively, we need to write $V=V_1V_2$, where $V_1$ decays like $\la x\ra^{-1-\eps}$, whereas the other should decay like $\la x\ra^{-\half-\eps}$.
Thus, in this case $\beta>\frac32$
is enough. On the other hand, in~\eqref{eq:one_der} one derivative may fall on
one of the $G$-terms at the ends. Then $V$ has to compensate for a $\frac32+\eps$ power
because of~\eqref{eq:Gest1}, and also a $\frac12+\eps$ power from
the limiting absorption principle.
Similarly with the other terms.

As far as $I^{+}(t,x,y)$ is concerned, note that on the support
of~$a^{\pm}_{x,y}(\lambda)$ the phase
$t\lambda^2+\lambda(|x|+|y|)$ has no critical point.  A single
integration by parts yields the bound $|I^{+}(t,x,y)|\les t^{-1}$
uniformly in~$x,y\in\R^2$, see~\eqref{eq:adec}.

In the case of $I^{-}(t,x,y)$ the phase $t\lambda^2-\lambda(|x|+|y|)$
has a unique critical point at $\lambda_0=(|x|+|y|)/(2t)$. If $\lambda_0\ll \lambda_1$, then
a single integration by parts again yields the bound of~$t^{-1}$.
If $\lambda_0\gtrsim \lambda_1$ then the bound
$\max(|x|,|y|)\gtrsim t$ is also true, and stationary phase contributes
$t^{-\half}(\la x\ra \la y\ra)^{-\half}\les t^{-1}$,
as desired.  To make this rigorous, apply Lemma~\ref{lem:sp}:
\begin{align*}
&|I^{-}(t,x,y)| \\
&\les |t|^{-1}\int \frac{|a^{-}_{x,y}(\lambda)|}{\delta^2+|\lambda-\lambda_0|^2} \,d\lambda
+ |t|^{-1}\int_{[|\lambda-\lambda_0|>\delta]} \frac{|\partial_\lambda a^{-}_{x,y}(\lambda)|}{|\lambda-\lambda_0|} \,d\lambda \\
&\les |t|^{-1} \int_0^\infty \frac{(1+\lambda)^{-2}(\la x\ra\la y\ra)^{-\half}}{\delta^2+|\lambda-\lambda_0|^2} \,d\lambda
+ |t|^{-1}\int_{[|\lambda-\lambda_0|>\delta]} \frac{(1+\lambda)^{-2}(\la x\ra\la y\ra)^{-\eps} }{|\lambda-\lambda_0|} \,d\lambda  \\
& \les |t|^{-1}
\end{align*}
since $\la x\ra+\la y\ra\gtrsim t$.
Note that when $0<t<1$ one has the better bound
$|I^{\pm}(t,x,y)| \les 1$ by~\eqref{eq:adec}.
\end{proof}

\section{Energies close to zero}
\label{sec:low}

The following lemma is a variant of the standard asymptotic expansion around zero energy of the
free resolvent on $\R^2$. The estimates on the error terms are written in a somewhat unusal form,
which is the one needed later in the proof.

\begin{lemma}
\label{lem:free_exp}
Let $R_0^{\pm}(\lambda^2)=(-\Laplace-(\lambda^2\pm i0))^{-1}$ be the free resolvent
in $\R^2$. Then, for all $\lambda>0$,
\begin{equation}
\label{eq:zero_exp}
R_0^{\pm}(\lambda^2)  = \Big[\pm\frac{i}{4} - \frac{1}{2\pi}\gamma - \frac{1}{2\pi} \log(\lambda/2)\Big]P_0
+ G_0  + E_0^{\pm}(\lambda).
\end{equation}
Here $P_0f:=\int_{\R^2} f(x)\,dx$, $G_0f(x)=-\frac{1}{2\pi}\int_{\R^2} \log|x-y|\,f(y)\,dy$, and the
error $E_0^{\pm}(\lambda)$ satisfies
\beeq
\label{eq:err}
\big\|\sup_{0<\lambda}\lambda^{-\half}|E_0^{\pm}(\lambda)(\cdot,\cdot)|\,\big \|+
\big\|\sup_{0<\lambda}\lambda^{\half} |\partial_\lambda E_0^{\pm}(\lambda)(\cdot,\cdot)|\,\big\|
\les 1
\eneq
with respect to the Hilbert-Schmidt
norm in $\bdd(L^{2,s}(\R^2),L^{2,-s}(\R^2))$ with $s>\frac32$.
\end{lemma}
\begin{proof}
One has, with $\lambda>0$,
\beeq
\label{eq:H0}
 R_0^{\pm}(\lambda^2)(x,y)=\pm \frac{i}{4} H_0^{\pm}(\lambda|x-y|)
\eneq
where the Hankel functions $H_0^{\pm}$ are
\begin{align*}
H_0^{\pm}(z) &= J_0(z)\pm iY_0(z) \\
& = 1 \pm i\frac{2}{\pi}\gamma \pm i\frac{2}{\pi} \log(z/2) + O(z^2\log z),\\
 \frac{d}{dz}H_0{^{\pm}}(z) &= \pm i\frac{2}{\pi z} +  O(z\log z).
\end{align*}
This is an expansion around $z=0$. Around $z=\infty$ the expansion is given by
\[ H_0^{\pm}(z)= \sqrt{\frac{2}{\pi z}}(a(z)\pm i b(z))e^{\pm i(z-\pi/4)},\]
with $a(z)=1-\frac{\alpha}{z^2}\pm\ldots$ and $b(z)=\frac{\beta}{z}\pm\ldots$.
Now let
\[
E_0^{\pm}(\lambda)(x,y):= R_0^{\pm}(\lambda^2)(x,y)- \Big[\pm\frac{i}{4} -
\frac{1}{2\pi}\gamma - \frac{1}{2\pi} \log(\lambda/2)\Big]+\frac{1}{2\pi}\log|x-y|.
\]
Then
\begin{align*}
&|E_0^{\pm}(\lambda)(x,y)| \\
 &\les \lambda^2|x-y|^2|\log(\lambda|x-y|)|\chi_{[\lambda|x-y|\le 1]}+
[1+\log (\lambda |x-y|)]\chi_{[\lambda|x-y|>1]} \\
&\les \lambda^{2\eps}|x-y|^{2\eps} |\log(\lambda|x-y|)|\chi_{[\lambda|x-y|\le 1]}+
[1+\log (\lambda |x-y|)]\chi_{[\lambda|x-y|>1]}.
\end{align*}
Hence
\[ \sup_{0<\lambda} \lambda^{-\eps} |E_0^{\pm}(\lambda)(x,y)| \les |x-y|^{\eps}.\]
Since the right-hand side has finite Hilbert-Schmidt norm as an operator
$L^{2,s}(\R^2) \to L^{2,-s}(\R^2))$ with $s>1+\eps$, we obtain the first part of~\eqref{eq:err}.
On the other hand,
\begin{align*}
& \lambda^{1-\eps}|\partial_\lambda E_0^{\pm}(\lambda)(x,y)| \\
 &\les \lambda^{2-\eps}|x-y|^2|\log(\lambda|x-y|)|\chi_{[\lambda|x-y|\le 1]} + [\lambda^{\half-\eps}|x-y|^{\half}+\lambda^{-\eps}]\chi_{[\lambda|x-y|>1]} \\
&\les |x-y|^{\eps}+\lambda^{\half-\eps} |x-y|^{\half}\chi_{[\lambda|x-y|>1]},
\end{align*}
and therefore, setting $\eps=\half$,
\[ \sup_{0<\lambda} \lambda^{\half} |\partial_\lambda E_0^{\pm}(\lambda)(x,y)| \les |x-y|^{\half}.\]
Since the right-hand side has finite Hilbert-Schmidt norm as an operator
$L^{2,s}(\R^2) \to L^{2,-s}(\R^2)$ with $s>\frac32$,
the lemma follows.
\end{proof}

\noindent Now let $V:\R^2\to\R$, $V\not\equiv0$,  satisfy $|V(x)|\les (1+|x|)^{-2\beta}$
for $\beta>\frac32$
(this condition arises because of the condition $s>\frac32$ in the previous lemma).
Following Jensen and Nenciu~\cite{JenNen} we set $U(x)=1$ if $V(x)\ge0$ and $U(x)=-1$ if $V(x)<0$.
Also, $v(x):=|V(x)|^{\half}$ decays like $(1+|x|)^{-\beta}$. The following corollary is therefore
an immediate consequence of Lemma~\ref{lem:free_exp}.

\begin{cor}
\label{cor:M}
For $\lambda>0$ define $M^{\pm}(\lambda):=U+vR_0^{\pm}(\lambda^2)v$. Let $P=\frac{v\la\cdot,v\ra}{\|V\|_1}$ denote
the orthogonal projection onto~$v$.
Then
\beeq
\label{eq:Mexp}
M^{\pm}(\lambda) = g^{\pm}(\lambda)P+U+vG_0v + vE_0^\pm(\lambda) v.
\eneq
Here $G_0,E_0^{\pm}(\lambda)$ are as in Lemma~\ref{lem:free_exp} and
$g^{\pm}(\lambda)=\|V\|_1\Big(\pm\frac{i}{4}- \frac{1}{2\pi}\gamma - \frac{1}{2\pi}\log(\lambda/2)\Big)$.
The remainders satisfy
\beeq
\label{eq:vEv}
 \| v\,\sup_{0<\lambda}\lambda^{-\half}|E_0^{\pm}(\lambda)(\cdot,\cdot)|\, v\|_{HS} +
\| v \sup_{0<\lambda}\lambda^{\half} |\partial_\lambda E_0^{\pm}(\lambda)(\cdot,\cdot)|\,  v\|_{HS} \les 1
\eneq
 with respect to the Hilbert-Schmidt norm on~$L^2(\R^2)$.
\end{cor}

\noindent The following definition is motivated by~\cite{JenNen}, cf.~the case of $S_1=0$
in their Theorem~6.2.

\begin{defi}
\label{defi:zero}  Let $Q=1-P$.
We say that zero is a {\em regular point} of the spectrum of $H=-\Laplace+V$ provided
$Q(U+vG_0v)Q$ is invertible on~$QL^2(\R^2)$. In that case set $D_0:=[Q(U+vG_0v)Q]^{-1}$
as an operator on~$QL^2(\R^2)$.
\end{defi}

\noindent Jensen and Nenciu show that $Q(U+vG_0v)\Phi=0$ with $\Phi\in QL^2(\R^2)$
implies that $\Phi=Uv\Psi$ where $H\Psi=0$ in the sense of distributions and $\Psi\in L^\infty(\R^2)$.
Thus Definition~\ref{defi:zero} captures what is sometimes described as absence of zero-energy
eigenfunctions and resonances.

\noindent
The following lemma is a technical statement that will be used repeatedly in our argument.

\begin{lemma}
\label{lem:absL2}
Let $D_0$ be as in Definition~\ref{defi:zero}.
Let $K$ be the kernel of the operator $QD_0Q$. Then the operator with kernel $|K|$ is
again $L^2$-bounded.
\end{lemma}
\begin{proof} For the purposes of this proof we introduce the following
terminology: A bounded operator $T$ on $L^2(\R^2)$ is called {\em absolutely bounded} if
the absolute value of its kernel gives rise to a bounded operator on $L^2(\R^2)$. Note that
Hilbert-Schmidt operators are absolutely bounded.

\noindent
Suppose $f\in QL^2(\R^2)$ such that $QUf=0$, $f\ne0$.
Then $Uf=cv$ for some scalar~$c\ne0$. Hence $f=cUv$ and $Pf=0$ requires
that $\la f,v\ra=c\la Uv,v\ra=c\int_{\R^2}V(x)\,dx=0$. Since this argument can be
reversed,
\[ \ker_{QL^2}(QUQ)=\{0\} \text{\ \ iff\ \ }\int_{\R^2} V(x)\,dx \ne 0.\]

\noindent {\em Case 1:} $\int_{\R^2} V(x)\,dx \ne0$.

\medskip
In this case we claim that $QUQ:QL^2(\R^2)\to QL^2(\R^2)$ is invertible. More precisely,
one checks that for any $g\in L^2$ with $Qg=g$
\[ f=Ug+ c_0 Uv \text{\ \ with\ \ }c_0=-\frac{\la Ug,v\ra}{\int_{\R^2} V(x)\,dx}\]
solves $QUQf=g$, $Qf=f$. It is evident from this explicit formula that $Q(QUQ)^{-1}Q$
is absolutely bounded. Moreover, on $QL^2$,
\beeq
\label{eq:rep1}
 [Q(U+vG_0v)Q]^{-1}=(QUQ)^{-1}[Q+QvG_0v(QUQ)^{-1}Q]^{-1}.
\eneq
Now $vG_0v$ is a Hilbert-Schmidt operator since $v$ decays faster than $(1+|x|)^{-1-\eps}$.
Hence $W:=QvG_0v(QUQ)^{-1}Q$ is also Hilbert-Schmidt.
Finally, as an identity on $QL^2$,
\[ [Q+QvG_0v(QUQ)^{-1}Q]^{-1}-Q=-[Q+W]^{-1}W. \]
Since the right-hand side is Hilbert-Schmidt, we see from~\eqref{eq:rep1}
that $[Q(U+vG_0v)Q]^{-1}$ is the composition of an absolutely bounded operator
with the sum of an absolutely bounded operator and a Hilbert-Schmidt operator.
Hence it is itself absolutely bounded, as claimed.

\medskip
\noindent {\em Case 2:} $\int_{\R^2} V(x)\,dx =0$.

\medskip In this case we remark that $0$ is an isolated point of the spectrum
of $QUQ$. Let $\pi_0$ denote the Riesz projection onto $\ker(QUQ)$ in $QL^2$.
From the preceding, $\pi_0(f)=\|V\|_1^{-1} \la f,Uv\ra Uv$.
Then $QUQ+\pi_0$ is invertible on~$QL^2(\R^2)$. In fact, one checks that an explicit solution of
\[ (QUQ+\pi_0)f=g \text{\ \ where\ \ } Qg=g,\; Qf=f\]
is given by
\[ f= Ug+c_1 v -c_1 Uv\text{\ \ with\ \ }c_1=-\frac{\la g,Uv\ra}{\int_{\R^2}|V(x)|\,dx}.\]
In view of this explicit expression, $(QUQ+\pi_0)^{-1}$ is absolutely bounded on $QL^2$.
Finally,  the identity
\[ [Q(U+vG_0v)Q]^{-1} = [QUQ+\pi_0+QvG_0vQ-\pi_0]^{-1} \]
on $QL^2$ allows one to repeat the same argument as in Case~1 and the lemma follows.
\end{proof}

The main technical result in Jensen and Nenciu~\cite{JenNen} is a formula
for the inverse of~$M^{\pm}(\lambda)^{-1}$. In the general case, this is complicated,
see their Theorem~6.2. But since we are imposing the condition of Definition~\ref{defi:zero},
it is relatively simple to compute that inverse,
see the following lemma.
Since we need somewhat stronger bounds on the error than those obtained in~\cite{JenNen},
we give all details. In particular, the proof requires Lemma~\ref{lem:absL2}.

\begin{lemma}
\label{lem:inv_exp}
Suppose that zero is a regular point of the spectrum of $H=-\Laplace+V$. Then
for some sufficiently small $\lambda_1>0$, the operators $M^{\pm}(\lambda)$ are invertible
for all $0<\lambda<\lambda_1$ as bounded operators on $L^2(\R^2)$, and one has the expansion
\beeq
\label{eq:M_inv_exp}
M^{\pm}(\lambda)^{-1} = h_{\pm}(\lambda)^{-1}S+QD_0Q+E^{\pm}(\lambda),
\eneq
where $h_{+}(\lambda)=a\log \lambda+z$, $a$ is real, $z$ complex, $a\ne0$, $\Im z\ne0$, and $h_{-}(\lambda)=\overline{h_+(\lambda)}$.
Moreover, $S$ is of finite rank and has a real-valued kernel,
and $E^{\pm}(\lambda)$ is a Hilbert-Schmidt operator
that satisfies the bound
\beeq
\label{eq:Eest}
\big \|\sup_{0<\lambda<\lambda_1} \lambda^{-\half} |E^{\pm}(\lambda)(\cdot,\cdot)|\,\big  \|_{HS}
+  \big \|\sup_{0<\lambda<\lambda_1} \lambda^{\half} |\partial_\lambda E^{\pm}(\lambda)(\cdot,\cdot)|\, \big \|_{HS} \les 1
\eneq
where the norm refers to the Hilbert-Schmidt norm on $L^2(\R^2)$. Finally, let $R_V^{\pm}(\lambda^2) = (-\Laplace+V-(\lambda^2\pm i0))^{-1}$.  Then
\beeq
R_V^{\pm}(\lambda^2)  = R_0^{\pm}(\lambda^2) -
R_0^{\pm}(\lambda^2)v M^{\pm}(\lambda)^{-1} v R_0^{\pm}(\lambda^2).
\label{eq:RV}
\eneq
This is to be understood as an identity between operators $L^{2,\half+\eps}(\R^2)\to L^{2,-\half-\eps}(\R^2)$ for some sufficiently small $\eps>0$.
\end{lemma}
\begin{proof}
For the purposes of this proof set $T=U+vG_0v$. By assumption, $QTQ$ is invertible on $QL^2(\R^2)$.
Moreover, by Corollary~\ref{cor:M}, with respect to the decomposition
$L^2(\R^2)=PL^2(\R^2)\oplus QL^2(\R^2)$,
\[
M^{\pm}(\lambda)= \left[ \begin{matrix} g^{\pm}(\lambda)P+PTP & PTQ \\
QTP & QTQ \end{matrix} \right] + vE_0^{\pm}(\lambda)v.
\]
Denote the matrix on the right-hand side
by~$A(\lambda)=\left[\begin{matrix} a_{11} & a_{12} \\ a_{21} &
a_{22} \end{matrix}\right]$. To invert $M^{\pm}(\lambda)$ and
thus~$A(\lambda)$, we use the well-known Fehsbach formula, see
eg.~Lemma~2.3 in~\cite{JenNen}. This requires that
$a:=(a_{11}-a_{12}a_{22}^{-1}a_{21})^{-1}$ exists, and in that
case \beeq \label{eq:Ainv}
 A(\lambda)^{-1} = \left[\begin{matrix} a & -aa_{12}a_{22}^{-1} \\ -a_{22}^{-1}a_{21}a &
a_{22}^{-1}a_{21}aa_{12}a_{22}^{-1}+a_{22}^{-1} \end{matrix}\right].
\eneq
Note that in our case, as an operator on the line $\Ran(P)=\{cv\::\:c\in\Compl\}$,
\[ a=h_{\pm}(\lambda)^{-1}P\text{\ \ where\ \ }h_{\pm}(\lambda):=g^{\pm}(\lambda)+\trace(PTP-PTQD_0QTP).\]
The trace is real-valued since $v$ is real-valued and since the kernel of $T$ is
real-valued. In view of the definition of $g^{\pm}(\lambda)$, $h_{\pm}(\lambda)\ne0$  provided $\lambda>0$
is sufficiently small. Moreover, by~\eqref{eq:Ainv} we see that
\[ A(\lambda)^{-1} = h_{\pm}(\lambda)^{-1}S+QD_0Q\]
where $S$ is of finite rank (in fact, the rank is at most two).
By the definition of $h_{\pm}(\lambda)$ and by Lemma~\ref{lem:absL2},
\beeq
\label{eq:Aest}
|A(\lambda)^{-1}(\cdot,\cdot)|+\lambda |\partial_\lambda A(\lambda)^{-1}(\cdot,\cdot)|
\les |S(\cdot,\cdot)|+ |QD_0Q(\cdot,\cdot)|,
\eneq
where the right-hand side is an $L^2$-bounded operator.
Now
\[ M^{\pm}(\lambda)^{-1} = A(\lambda)^{-1}[1+vE_0^{\pm}(\lambda)v A(\lambda)^{-1}]^{-1}.\]
The second inverse on the right-hand side exists for small $\lambda$ since then
\[ \|vE_0^{\pm}(\lambda)v A(\lambda)^{-1}\| < \half,\] see~\eqref{eq:vEv}.
Moreover, writing out $E^{\pm}(\lambda)$ as a Neuman series, one obtains~\eqref{eq:Eest}
from~\eqref{eq:vEv} and~\eqref{eq:Aest} by termwise estimation.

\noindent
Finally, \eqref{eq:RV} is the well-known symmetric resolvent expansion which follows
easily from
\begin{align*}
& (I-Uv(-\Laplace+V-z)^{-1}v)(I+Uv(-\Laplace-z)^{-1} v) = I {\ \ \ or\ \ } \\
& V(-\Laplace+V-z)^{-1}V =V-v(U+v(-\Laplace-z)^{-1}v)^{-1}v
\end{align*}
for $\Im z>0$. Passing to the limit $\Im z\to 0$ now leads to~\eqref{eq:RV} via an
application of the resolvent identity and the limiting absorption principle, cf.~\eqref{eq:lim_ap}.
\end{proof}

\begin{cor}
\label{cor:RV_exp}
Let zero be a regular point of the spectrum of $H=-\Laplace+V$. Then
\begin{align}
R_V^{\pm}(\lambda^2) &= R_0^{\pm}(\lambda^2) - h_{\pm}(\lambda)^{-1}\,R_0^{\pm}(\lambda^2)v S v R_0^{\pm}(\lambda^2) \nn\\
& -R_0^{\pm}(\lambda^2)v QD_0Q v R_0^{\pm}(\lambda^2)- R_0^{\pm}(\lambda^2)vE^{\pm}(\lambda) v R_0^{\pm}(\lambda^2)
\label{eq:RV_exp}
\end{align}
where $S$ and $E^{\pm}(\lambda)$ are as in the previous lemma.
This is to be understood as an identity between operators $L^{2,\half+\eps}(\R^2)\to L^{2,-\half-\eps}(\R^2)$ for small $\eps>0$,
i.e., as in the limiting absorption principle~\eqref{eq:lim_ap}.
\end{cor}
\begin{proof} This is an immediate consequence of Corollary~\ref{cor:M} and Lemma~\ref{lem:inv_exp}.
\end{proof}

We now turn to decay estimates.

\begin{prop}
\label{prop:low_decay}
Let $\chi$ be a smooth cut-off function on the line with $\chi(\lambda)=1$
if $\lambda\le\lambda_1$ and $\chi(\lambda)=0$ if $\lambda\ge 2\lambda_1$ where $\lambda_1>0$ is a small
constant. Assume that zero is a regular point of the spectrum of $H=-\Laplace+V$. Then
\bear
|\la e^{itH}\chi(\sqrt{H})P_{ac}(H) f,g \ra | &=& \frac{1}{\pi}\Bigl| \int_0^\infty e^{it\lambda^2}\lambda \chi(\lambda)\big\la[R_V^{+}(\lambda^2)-R_V^{-}(\lambda^2)]f,g\big\ra \, d\lambda \Bigr| \nn \\
& \le&  C |t|^{-1}\|f\|_1\|g\|_1 \label{eq:low_dec}
\eear
for every $f,g\in\calS(\R^2)$. Here $C$ is a constant that only depends on $V$ and $\chi$.
\end{prop}

The proof of Proposition~\ref{prop:low_decay} is based on the expansion of $R_V^{\pm}(\lambda^2)$
stated in Corollary~\ref{cor:RV_exp}. Each of the four terms on the right-hand side of~\eqref{eq:RV_exp}
requires a separate argument.
We begin with the free case.

\begin{lemma}
\label{lem:I}
$H_0=-\Laplace$ satisfies
\[ |\la e^{itH_0}\chi(\sqrt{H_0})P_{ac}(H) f,g \ra | \le C|t|^{-1}\|f\|_1\|g\|_1 \]
for all $f,g\in\calS(\R^2)$.
\end{lemma}
\begin{proof}
This follows immediately from the standard bound $\|e^{itH_0}f\|_\infty \les |t|^{-1} \|f\|_1$
and the fact that $\chi(\sqrt{H_0})$ and $\Pac(H)$ are bounded on $L^1(\R^2)$ (for the latter,
use that the number of negative bound states is finite~\cite{mihai}, as well as
that the eigenfunctions are exponentially decaying by Agmon's bound, and therefore in $L^1(\R^2)$.
Moreover, they are in $L^\infty(\R^2)$ by Sobolev imbedding).
Alternatively, one can give a self-contained proof via stationary phase. Indeed, from~\eqref{eq:H0}
\[ R_0^{+}(\lambda^2)(x,y)-R_0^{-}(\lambda^2)(x,y)=\frac{i}{2}J_0(\lambda|x-y|).\]
Thus
\begin{align*}
& |\la e^{itH_0}\chi(\sqrt{H_0})P_{ac}(H) f,g \ra | \\
& \le \int_{\R^4}\frac{1}{2\pi}
\Bigl| \int_0^\infty e^{it\lambda^2}\lambda \chi(\lambda) J_0(\lambda|x-y|) \, d\lambda \Bigr|
|P_{ac}(H) f(x)||g(y)|\,dxdy.
\end{align*}
Now $J_0(u)=e^{iu}\omega_+(u)+e^{-iu}\omega_-(u)$ where $|\omega_{\pm}(u)|\les (1+|u|)^{-\half}$.
Therefore,
\begin{align}
& \Bigl| \int_0^\infty e^{it\lambda^2}\lambda \chi(\lambda) J_0(\lambda|x-y|) \, d\lambda \Bigr| \nn \\
& \les \Bigl| \int_0^\infty e^{i[t\lambda^2-\lambda|x-y|]}\lambda \chi(\lambda) \omega_{+}(\lambda|x-y|) \, d\lambda \Bigr| \label{eq:om-} \\
& + \Bigl| \int_0^\infty e^{i[t\lambda^2+\lambda|x-y|]}\lambda \chi(\lambda) \omega_{-}(\lambda|x-y|) \, d\lambda \Bigr|. \label{eq:om+}
\end{align}
Let $t>0$.
The phase in~\eqref{eq:om-} has a stationary point $\lambda_0=\frac{|x-y|}{2t}$. Hence that integral
is
\[ \les t^{-\half}\lambda_0 (1+\lambda_0^2 t)^{-\half} \les t^{-1}\]
by stationary phase (we leave it to the reader to fill in the remaining details here).
The integral in~\eqref{eq:om+} can be estimated directly by means of  integration by parts.
\end{proof}

The following lemmas deal with the contribution of the term containing $QD_0Q$ in~\eqref{eq:RV_exp}.
In what follows it will be assumed that zero is a regular point of the spectrum of $H=-\Laplace+V$.

\begin{lemma}
\label{lem:II_1} Let $(QD_0Q)(\cdot,\cdot)$ denote the kernel of $QD_0Q$. There is the bound
\bear
&& \Bigl|\int_{\R^8}\int_0^\infty e^{it\lambda^2} \lambda \chi(\lambda) \chi(\lambda|x-x_1|) Y_0(\lambda|x-x_1|)
v(x_1)(QD_0Q)(x_1,y_1) v(y_1) \nn \\
&& J_0(\lambda|y_1-y|)\chi(\lambda|y_1-y|)\,d\lambda\,f(x)g(y)\,dx_1dy_1\,dxdy\Bigr|
\le C\, |t|^{-1}\|f\|_1\|g\|_1 \label{eq:QD0Q_one}
\eear
with a constant that only depends on $V$.
\end{lemma}
\begin{proof}
 We  make the following claim:
\bear
&& \Bigl|\int_0^\infty e^{it\lambda^2} \lambda \chi(\lambda) \Big[\chi(\lambda|x-x_1|) Y_0(\lambda|x-x_1|)
-\frac{2}{\pi}\chi(\lambda(1+|x|))\log(\lambda(1+|x|))\Big] \nn \\
&&  J_0(\lambda|y_1-y|)\chi(\lambda|y_1-y|)\,d\lambda\,\Bigr|
\le C\, |t|^{-1}(1+\log^+|x_1|+\log^{-}|x-x_1|)
\label{eq:claim1}
\eear
for all $x,x_1,y,y_1\in\R^2$.
Let
\beeq
\label{eq:kdef}
k(x,x_1):=1+\log^+|x_1|+\log^{-}|x-x_1|.
\eneq
If~\eqref{eq:claim1} holds,
then the left-hand side of~\eqref{eq:QD0Q_one} is
\begin{align*}
&\les |t|^{-1}\int_{\R^8} k(x,x_1)v(x_1)|(QD_0Q)(x_1,y_1)| v(y_1)|f(x)||g(y)|\,dx_1dy_1dxdy \\
&\les |t|^{-1} \sup_{x\in\R^2} \Big(\int_{\R^2} k(x,x_1)^2|V|(x_1)\,dx_1\Bigr)^{\half} \|\,|QD_0Q|\,\|_{2\to2} \|V\|_1^{\half}\|f\|_1\|g\|_1 \\
&\les |t|^{-1} \|f\|_1\|g\|_1,
\end{align*}
as desired. To see this, observe firstly that
\beeq
\label{eq:0int}
 \int_{\R^4} v(x)(QD_0Q)(x,y)h(y)\,dxdy =0
\eneq
for any $h\in L^2(\R^2)$. Secondly,
 use Lemma~\ref{lem:absL2} to control the $L^2$-operator norm of the kernel $|QD_0Q|$.
To prove~\eqref{eq:claim1}, let
\begin{align}
 F(\lambda,x,x_1) &:=\chi(\lambda|x-x_1|) Y_0(\lambda|x-x_1|)
-\frac{2}{\pi}\chi(\lambda(1+|x|))\log(\lambda(1+|x|)) \label{eq:Fdef}\\
 G(\lambda,y_1,y) &:= J_0(\lambda|y_1-y|)\chi(\lambda|y_1-y|).\label{eq:Gdef}
\end{align}
If we choose $1>\lambda_1>0$ so that $2\lambda_1$ lies to the left of the first
zero of~$J_0$, then $G(\lambda,y_1,y)$ is nonincreasing in~$\lambda$ (recall the definition
of $\chi$ in Proposition~\ref{prop:low_decay}). Moreover, in
that case $0\le G\le 1$ for all choices of arguments.
Recall that $J_0(z)=1+O(z^2)$ and
\begin{align}
Y_0(z) &=\frac{2}{\pi}(\log z + c)J_0(z)+r(z) \label{eq:Y0} \\
Y_0'(z) &= \frac{2}{\pi z}J_0(z)+\frac{2}{\pi}(\log z+c)J_0'(z)+r'(z)=\frac{2}{\pi z}+g(z) \label{eq:Y0'}
\end{align}
where $r(z)$ is analytic for all $z$  and $g(z)$ bounded on $(0,\infty)$, say.
Hence one has $F(0+,x,x_1)=\frac{2}{\pi}c+\frac{2}{\pi}\log\frac{|x-x_1|}{1+|x|}$, and $G(0,y_1,y)=1$.
It is easy to check that
\beeq
\label{eq:logbd}
\Bigl| \log\frac{|x-x_1|}{1+|x|} \Bigr| \les 1+\log^+ |x_1| + \log^{-}|x-x_1|=k(x,x_1).
\eneq
Indeed, if $|x|\ge 2|x_1|$, then
\[ \frac14\chi_{[|x|\ge1]} + \half|x-x_1|\chi_{[|x|\le1]} \le \frac{|x-x_1|}{1+|x|} \le \frac{2|x|}{1+|x|} \le 2.\]
On the other hand, if $|x|<2|x_1|$, then
\[ \frac{\min(1,|x-x_1|)}{1+2|x_1|}\le \frac{|x-x_1|}{1+|x|} \le \frac{3|x_1|}{1+|x|}\le 3|x_1|,\]
and~\eqref{eq:logbd} follows. Integrating by parts inside the integral in~\eqref{eq:claim1}
therefore leads to the estimate
\bear
\eqref{eq:claim1} &\les& |t|^{-1} k(x,x_1) + |t|^{-1}\int_0^\infty  |\chi'(\lambda)| |F(\lambda,x,x_1)|
|G(\lambda,y_1,y)| \,d\lambda \nn \\
&& + |t|^{-1}\int_0^\infty   |\partial_\lambda F(\lambda,x,x_1)|
|G(\lambda,y_1,y)| \,d\lambda \label{eq:DF} \\
&& + |t|^{-1}\int_0^\infty   |F(\lambda,x,x_1)|
|\partial_\lambda G(\lambda,y_1,y)| \,d\lambda. \label{eq:DG}
\eear
Recall that the support of $\chi'$ is contained inside $[\lambda_1,2\lambda_1]$. Thus
the integral involving $\chi'(\lambda)$ is easily seen to be
\[ \les \sup_{\lambda\sim \lambda_1} |F(\lambda,x,x_1)| |G(\lambda,y_1,y)| \les 1+\log^{-}|x-x_1|,\]
cf.~\eqref{eq:Y0}.
With the notation of~\eqref{eq:Y0'},
\begin{align}
& \partial_\lambda F(\lambda,x,x_1) = \frac{2}{\pi}\frac{1}{\lambda} \Bigl[\chi(\lambda|x-x_1|)-
\chi(\lambda(1+|x|))\Bigr] \nn \\
& + |x-x_1|\chi'(\lambda|x-x_1|) Y_0(\lambda|x-x_1|)
+ |x-x_1| \chi(\lambda|x-x_1|)g(\lambda|x-x_1|) \nn \\
& - \frac{2}{\pi}\chi'(\lambda(1+|x|))(1+|x|)\log(\lambda(1+|x|)). \label{eq:F'}
\end{align}
Hence,
\begin{align*}
\eqref{eq:DF} &\les |t|^{-1}\int_0^\infty  |\chi(\lambda|x-x_1|)-
\chi(\lambda(1+|x|))| \lambda^{-1}\, d\lambda  \\
& +  |t|^{-1}\int_0^\infty |x-x_1|[|\chi'(\lambda|x-x_1|)|+ \chi(\lambda|x-x_1|) ] \,d\lambda \\
& +  |t|^{-1}\int_0^\infty (1+|x|)|\chi'(\lambda(1+|x|))|\, d\lambda \\
&\les  |t|^{-1}\Big[ 1 + \log^+ \Big(\frac{2|x-x_1|}{1+|x|}\Big) + \log^+ \Big(\frac{2(1+|x|)}{|x-x_1|}\Big)\Big]
\les |t|^{-1}\,k(x,x_1),
\end{align*}
where we used \eqref{eq:logbd} in the last step.
In passing, we note that we have shown the following:
\beeq
\label{eq:Fbound}
\sup_{0\le \lambda\le 1} |F(\lambda,x,x_1)| \le |F(0,x,x_1)|+\int_0^1 |\partial_\lambda F(\lambda,x,x_1)|\,d\lambda \les k(x,x_1).
\eneq
As observed previously, $\partial_\lambda G$ has a definite sign. Moreover,
$F(\lambda,x,x_1)$ only has a finite number of zeros in~$\lambda$. Hence, one
can break up the integral~\eqref{eq:DG} into finitely many disjoint intervals,
remove the absolute values on each of them, and then integrate by parts. The only
boundary contribution occurs at $\lambda=0$, for which we have already obtained the desired bound.
Otherwise, the remaining integral is bounded above by~\eqref{eq:DF}, and we are done.
\end{proof}

The following lemma deals with an integral very much like the one
in~\eqref{eq:QD0Q_one}. The difference here is that we consider
the contribution from large arguments inside $J_0$, which makes it
necessary to exploit the oscillations of $J_0$. This will be done
by means of Lemma~\ref{lem:sp}.

\begin{lemma}
\label{lem:II_2} Let $(QD_0Q)(\cdot,\cdot)$ denote the kernel of $QD_0Q$. Let $\tchi=1-\chi$. Then there is the bound
\begin{align}
& \Bigl|\int_{\R^8}\int_0^\infty e^{it\lambda^2} \lambda \chi(\lambda) \chi(\lambda|x-x_1|) Y_0(\lambda|x-x_1|)
v(x_1)(QD_0Q)(x_1,y_1) v(y_1) \nn \\
& J_0(\lambda|y_1-y|)\tchi(\lambda|y_1-y|)\,d\lambda\,f(x)g(y)\,dx_1dy_1\,dxdy\Bigr|
\le C\, |t|^{-1}\|f\|_1\|g\|_1 \label{eq:QD0Q_two}
\end{align}
with a constant that only depends on $V$. The same statement holds with the role of the cut-offs
interchanged, i.e., with $\tchi(\lambda|x-x_1|)$
and $\chi(\lambda|y-y_1|)$.
\end{lemma}
\begin{proof} As usual,
\beeq
\label{eq:J0_decomp}
 J_0(y)=e^{iy}\omega_+(y) + e^{-iy}\omega_{-}(y)
\eneq
where $|\omega^{(\ell)}_{\pm}(y)|\les (1+|y|)^{-\half-\ell}$ for all $\ell\ge0$.
Correspondingly, there will be two contributions to~\eqref{eq:QD0Q_two}. We start with the phase
$\phi_{-}(\lambda)=\lambda^2-\lambda|y-y_1|t^{-1}$ which has a  critical point at $\lambda_0=\frac{|y-y_1|}{2t}$.
In that case we claim that
\begin{align}
& \Bigl|\int_0^\infty e^{it\phi_{-}(\lambda)} \lambda \chi(\lambda) \Big[\chi(\lambda|x-x_1|) Y_0(\lambda|x-x_1|)
-\frac{2}{\pi}\chi(\lambda(1+|x|))\log(\lambda(1+|x|))\Big] \nn \\
&  \omega_{-}(\lambda|y_1-y|)\tchi(\lambda|y_1-y|)\,d\lambda\,\Bigr|
\le C\, |t|^{-1} k(x,x_1)
\label{eq:claim2}
\end{align}
for all $x,x_1,y,y_1\in\R^2$. Here $k(x,x_1)$ is as in~\eqref{eq:kdef}. Moreover, as in the
previous proof, this bound will lead to the desired estimate in~\eqref{eq:QD0Q_two} in view
of~\eqref{eq:0int}.
With $F(\lambda,x,x_1)$ as in~\eqref{eq:Fdef}, set
\beeq
\label{eq:amin}
a(\lambda):=\lambda\chi(\lambda)\omega_{-}(\lambda|y-y_1|)\tchi(\lambda|y-y_1|)F(\lambda,x,x_1)
\eneq
where we suppress the other variables inside $a$. By Lemma~\ref{lem:sp},
\beeq
 \Big|\int_0^\infty e^{it\phi_{-}(\lambda)} a(\lambda)\,d\lambda\Big|
 \les |t|^{-1} \int_{-\infty}^\infty \Biggl( \frac{|a(\lambda)|}{\delta^2+|\lambda-\lambda_0|^2}
+\frac{|a'(\lambda)|}{|\lambda-\lambda_0|}\chi_{[|\lambda-\lambda_0|>\delta]} \Biggr) \, d\lambda.
\label{eq:2a}
\eneq
To establish our claim we therefore need to show that the integral in~\eqref{eq:2a} is
$\les k(x,x_1)$. Using~\eqref{eq:Fbound} one concludes
\begin{align}
& \int_{-\infty}^\infty  \frac{|a(\lambda)|}{\delta^2+|\lambda-\lambda_0|^2} \, d\lambda \nn \\
&\les k(x,x_1) \int \frac{\lambda \chi(\lambda) |\omega_{-}(\lambda|y-y_1|)|\tchi(\lambda|y-y_1|)}{\delta^2+|\lambda-\lambda_0|^2} \, d\lambda \nn  \\
& \les k(x,x_1)|y-y_1|^{-\half}\int_{c|y-y_1|^{-1}}^1 \frac{\sqrt{\lambda}}{\delta^2+|\lambda-\lambda_0|^2} \, d\lambda. \label{eq:reduc1}
\end{align}
Now suppose first that $\lambda_0\gtrsim \delta$, which is the same as $|y-y_1|\delta\gtrsim 1$.
Then
\begin{align*}
& |y-y_1|^{-\half}\int_{c|y-y_1|^{-1}}^1 \frac{\sqrt{\lambda}}{\delta^2+|\lambda-\lambda_0|^2} \, d\lambda\\
&\les |y-y_1|^{-\half}\Bigl\{\int_0^1 \frac{\sqrt{\lambda_0}}{\delta^2+|\lambda-\lambda_0|^2} \, d\lambda
+\int_0^1 \frac{\sqrt{|\lambda-\lambda_0|}}{\delta^2+|\lambda-\lambda_0|^2} \, d\lambda \Bigr\} \\
&\les |y-y_1|^{-\half}\Bigl\{\sqrt{\lambda_0}\delta^{-1}+\delta^{-\half} \Bigr\}
 \les 1 + (|y-y_1|\delta)^{-1} \les 1,
\end{align*}
as desired. On the other hand, if $\lambda_0 \ll \delta$, then also $|y-y_1|\delta\ll 1$ and
thus
\begin{align*}
& |y-y_1|^{-\half}\int_{c|y-y_1|^{-1}}^1
\frac{\sqrt{\lambda}}{\delta^2+|\lambda-\lambda_0|^2} \, d\lambda\\
&\les |y-y_1|^{-\half} \int_{c|y-y_1|^{-1}}^1 \lambda^{-\frac32}\,d\lambda \les 1.
\end{align*}
It remains to bound the contribution of the term involving $a'(\lambda)$ in~\eqref{eq:2a}.
Inspection of~\eqref{eq:F'} reveals that $|\partial_\lambda F(\lambda,x,x_1)|\les \lambda^{-1}$.
Combining this with~\eqref{eq:Fbound} yields
\begin{align}
|a'(\lambda)| & \les  k(x,x_1)\Big[(\chi(\lambda)+\lambda|\chi'(\lambda)|)(\lambda|y-y_1|)^{-\half}\tchi(\lambda|y-y_1|) \nn \\
& \qquad + \chi(\lambda)|\chi'(\lambda|y-y_1|)| \Big]. \label{eq:a'}
\end{align}
We start with the second term in~\eqref{eq:a'}. Its contribution to the integral in~\eqref{eq:2a}
is
\beeq
\les \int_0^1 \chi_{[|\lambda-\lambda_0|>\delta]} \frac{|\chi'(\lambda|y-y_1|)|}{|\lambda-\lambda_0|}
\, d\lambda.
\label{eq:a'1}
\eneq
The integration region here is contained inside an interval of the form
$[c_1|y-y_1|^{-1},c_2|y-y_1|^{-1}]$ where $c_1,c_2$ are some positive constants.
If $\lambda_0\asymp |y-y_1|^{-1}$, then also $|y-y_1|\delta\asymp1$. Hence in this case
\[ \eqref{eq:a'1} \les \log (1+\delta^{-1}|y-y_1|^{-1}) \les 1.\]
If on the other hand either $\lambda_0\gg |y-y_1|^{-1}$, or $\lambda_0\ll |y-y_1|^{-1}$,
then
\[ \eqref{eq:a'1} \les \Big|\log \Big( \frac{c_2|y-y_1|^{-1}-\lambda_0}{c_1|y-y_1|^{-1}-\lambda_0}\Big)
\Big| \les 1.
\]
It remains to consider the first term in~\eqref{eq:a'}.
Its contribution to the integral in~\eqref{eq:2a} is
\beeq
\les  |y-y_1|^{-\half} \int_{[\lambda|y-y_1|\gtrsim 1]}\chi_{[|\lambda-\lambda_0|>\delta]}\,
 \frac{d\lambda}{|\lambda-\lambda_0|\sqrt{\lambda}}.
\label{eq:a'2}
\eneq
If $\lambda_0\ll |y-y_1|^{-1}$, then
\[ \eqref{eq:a'2} \les |y-y_1|^{-\half} \int_{[\lambda|y-y_1|\gtrsim 1]}\frac{d\lambda}{\lambda^{\frac32}}
\les 1.\]
If, on the other hand, $\lambda_0 \gtrsim |y-y_1|^{-1}$, then
\begin{align*}
\eqref{eq:a'2} &\les |y-y_1|^{-\half} \Bigl\{ \int_0^{\half\lambda_0} \frac{d\lambda}{\lambda_0\sqrt{\lambda}} + \int_{\half\lambda_0}^1 \chi_{[|\lambda-\lambda_0|>\delta]} \frac{d\lambda}{|\lambda-\lambda_0|^{\frac32}} \Bigr\} \\
&\les |y-y_1|^{-\half}(\lambda_0^{-\half}+\delta^{-\half}) \les 1,
\end{align*}
as desired. In the last line we used that $\lambda_0 \gtrsim |y-y_1|^{-1}$ is
the same as $|y-y_1|\delta \gtrsim 1$.
This concludes the proof of claim~\eqref{eq:claim2}.

It remains to consider the phase $\phi_+(\lambda)=\lambda^2+t^{-1}|y-y_1|\lambda$.
The corresponding estimate is
\begin{align}
& \Bigl|\int_0^\infty e^{it\phi_{+}(\lambda)} \lambda \chi(\lambda) \Big[\chi(\lambda|x-x_1|) Y_0(\lambda|x-x_1|)
-\frac{2}{\pi}\chi(\lambda(1+|x|))\log(\lambda(1+|x|))\Big] \nn \\
&  \omega_{+}(\lambda|y_1-y|)\tchi(\lambda|y_1-y|)\,d\lambda\,\Bigr|
\le C\, |t|^{-1} k(x,x_1)
\label{eq:claim3}
\end{align}
for all $x,x_1,y,y_1\in\R^2$. Setting
\beeq
\label{eq:aplus}
a(\lambda):=\lambda\chi(\lambda)\omega_{+}(\lambda|y-y_1)\tchi(\lambda|y-y_1|)F(\lambda,x,x_1)
\eneq
a single integration by parts in the left-hand side of~\eqref{eq:claim3} yields
\beeq
\label{eq:a'again}
\eqref{eq:claim3} \les |t|^{-1} \int_0^\infty \frac{|a(\lambda)|}{|\phi_+'(\lambda)|^2}\, d\lambda
+ |t|^{-1} \int_0^\infty \frac{|a'(\lambda)|}{|\phi_+'(\lambda)|}\, d\lambda.
\eneq
As before,  $\lambda_0=\frac{|y-y_1|}{2t}$. Then
\begin{align*}
\int_0^\infty \frac{|a(\lambda)|}{|\phi_+'(\lambda)|^2}\, d\lambda
&\les k(x,x_1)|y-y_1|^{-\half} \int_0^\infty \lambda^{\half}(\lambda^2+\lambda_0^2)^{-1}\chi_{[\lambda|y-y_1|\gtrsim1]}\, d\lambda \\
&\les
 k(x,x_1)\int_0^\infty \lambda^{-\frac32} |y-y_1|^{-\half} \chi_{[\lambda|y-y_1|\gtrsim1]}\, d\lambda \\
& \les  k(x,x_1).
\end{align*}
To estimate the second integral in~\eqref{eq:a'again}, we use~\eqref{eq:a'} which remains
valid with $\omega_+$. Hence
\begin{align*}
\int_0^\infty \frac{|a'(\lambda)|}{|\phi_+'(\lambda)|}\, d\lambda
&\les  k(x,x_1) |y-y_1|^{-\half} \int_0^\infty \lambda^{-\half}(\lambda+\lambda_0)^{-1}\chi_{[\lambda|y-y_1|\gtrsim1]}\, d\lambda \\
&\les
k(x,x_1)\int_0^\infty \lambda^{-\frac32} |y-y_1|^{-\half} \chi_{[\lambda|y-y_1|\gtrsim1]}\, d\lambda \\
& \les k(x,x_1).
\end{align*}
In view of the preceding, $\eqref{eq:a'again} \les |t|^{-1}\,k(x,x_1)$. Hence~\eqref{eq:claim3}
holds and~\eqref{eq:QD0Q_two} has been proved.

\noindent The final statement about interchanging
the roles of $\chi$ and $\tchi$ is implicit in the previous proof. Indeed, \eqref{eq:J0_decomp}
holds equally well for $Y_0$ instead of~$J_0$. Moreover, one  replaces
$F(\lambda,x,x_1)$ with $G(\lambda,y,y_1)$, see~\eqref{eq:Gdef}, and the bound~\eqref{eq:Fbound}
with the trivial one $0\le G\le 1$.
We skip the details.
\end{proof}

The final lemma dealing with $QD_0Q$ controls the contributions of those
$\lambda$ for which both resolvents on either side of $vQD_0Qv$ are evaluated
at arguments of size $\gtrsim 1$. In this case it will be convenient to work
with the full kernel of the resolvents, i.e., the Hankel functions without splitting
them into $J_0$ and $Y_0$.

\begin{lemma}
\label{lem:II_3} Let $(QD_0Q)(\cdot,\cdot)$ denote the kernel of $QD_0Q$ and set $\tchi=1-\chi$. There is the bound
\begin{align}
& \Bigl|\int_{\R^4}\int_0^\infty\int_{\R^4} e^{it\lambda^2} \lambda \chi(\lambda) \tchi(\lambda|x-x_1|) H^{\pm}_0(\lambda|x-x_1|)
v(x_1)(QD_0Q)(x_1,y_1) v(y_1) \nn \\
& H^{\pm}_0(\lambda|y_1-y|)\tchi(\lambda|y_1-y|)\,d\lambda\,f(x)g(y)\,dx_1dy_1\,dxdy\Bigr|
\le C\, |t|^{-1}\|f\|_1\|g\|_1 \label{eq:QD0Q_three}
\end{align}
with a constant that only depends on $V$.
\end{lemma}
\begin{proof}
One has
\beeq
\label{eq:ompm}
 H_0^+(y)\tchi(y)=e^{iy}\omega_+(y) \text{\ \ and\ \ } H_0^-(y)\tchi(y)=e^{-iy}\omega_{-}(y)
\eneq
where $\omega_{-}=\overline{\omega_{+}}$, and
$|\omega_{\pm}^{(\ell)}(y)|\les (1+|y|)^{-\half-\ell}$ for all $\ell\ge0$ (the reader should note that we are slightly abusing notation here, since $\omega_{\pm}$ already
appeared as the decay factors of $J_0$ -- but this abuse of notation is of no consequence).
Correspondingly, there will be two phases to consider in~\eqref{eq:QD0Q_three},
namely
\[ \phi_{\pm}(\lambda) = \lambda^2 \pm \lambda\frac{|x-x_1|+|y-y_1|}{t}. \]
Set $p=|x-x_1|$ and $q=|y-y_1|$ for simplicity. We may assume that $p>0$ and $q>0$.
We claim that
\beeq
\label{eq:claim4}
\Bigl|\int_0^\infty e^{it\phi_{\pm}(\lambda)} \lambda \chi(\lambda) \tchi(\lambda p) \omega_{\pm}(p \lambda)
\tchi(\lambda q) \omega_{\pm}(q \lambda)\,d\lambda \Bigr| \les |t|^{-1},
\eneq
uniformly in $p,q>0$. The phase $\phi_{-}$ has a critical point at
\[ \lambda_0 = \frac{p+q}{2t}. \]
Let $a_{\pm}(\lambda)=\lambda \chi(\lambda) \tchi(\lambda p) \omega_{\pm}(p \lambda)
\omega_{\pm}(q \lambda)\tchi(\lambda q)$. Then by Lemma~\ref{lem:sp},
\begin{align}
\Bigl|\int_0^\infty e^{it\phi_{-}(\lambda)} a_{-}(\lambda)\,d\lambda \Bigr|
& \les |t|^{-1}\int_0^\infty \frac{|a_{-}(\lambda)|}{\delta^2+|\lambda-\lambda_0|^2}\,d\lambda \nn \\
& \qquad +|t|^{-1}\int_0^\infty \frac{|a_{-}'(\lambda)|}{|\lambda-\lambda_0|}\chi_{[|\lambda-\lambda_0|>\delta]}\,d\lambda. \label{eq:sp_1}
\end{align}
The integral involving $a_{-}(\lambda)$ is
\begin{align*}
& \les (pq)^{-\half} \int_{c(p^{-1}+q^{-1})}^1 \frac{d\lambda}{\delta^2+|\lambda-\lambda_0|^2} \\
& \les (pq)^{-\half} \Bigl(\delta^{-1}\chi_{[\lambda_0\gtrsim p^{-1}+q^{-1}]} + (p^{-1}+q^{-1})^{-1}\Bigr)
 \les 1.
\end{align*}
Here we used that $\lambda_0\gtrsim p^{-1}+q^{-1}$ is the same as $pq\gtrsim t$ or $pq\delta^2\gtrsim 1$,
as well as the bound
\[ (pq)^{-\half} (p^{-1}+q^{-1})^{-1} = \frac{\sqrt{pq}}{p+q} \lesssim 1. \]
Since
\[ |a_{-}'(\lambda)| \les (pq)^{-\half}\lambda^{-1}\chi_{[\lambda\gtrsim p^{-1}+q^{-1}]} \chi(\lambda),\]
the integral involving $a_{-}'(\lambda)$ in~\eqref{eq:sp_1} is
\beeq
\label{eq:a'3}
 \les (pq)^{-\half} \int_{c(p^{-1}+q^{-1})}^1 \chi_{[|\lambda-\lambda_0|>\delta]}\frac{d\lambda}{\lambda|\lambda-\lambda_0|} .
\eneq
Now suppose that $\lambda_0\ll \delta$. Then $|\lambda-\lambda_0|>\delta$ implies that
$\lambda-\lambda_0\gtrsim \lambda$. It follows that
\[ \eqref{eq:a'3}  \lesssim (pq)^{-\half} \int_{c(p^{-1}+q^{-1})}^1 \frac{d\lambda}{\lambda^2}
\lesssim (pq)^{-\half} (p^{-1}+q^{-1})^{-1} \lesssim 1.
\]
On the other hand, if $\lambda_0\gtrsim \delta$ which is the same as $(p+q)\delta\gtrsim1$,
then by Cauchy-Schwarz
\begin{align*}
 \eqref{eq:a'3}  &\lesssim (pq)^{-\half} \Big(\int_{c(p^{-1}+q^{-1})}^1 \frac{d\lambda}{\lambda^2}\Big)^{\half}  \Big(\int_{0}^1\chi_{[|\lambda-\lambda_0|>\delta]} \frac{d\lambda}{|\lambda-\lambda_0|^2}\Big)^{\half}  \\
& \lesssim (pq)^{-\half} (p^{-1}+q^{-1})^{-\half} \delta^{-\half} = (\delta(p+q))^{-\half} \lesssim 1.
\end{align*}
Hence \eqref{eq:claim4} holds for the phase $\phi_{-}$.

We now turn to $\phi_{+}$. By inspection,
\[
|a_+(\lambda)|\les (pq)^{-\half}\chi_{[\lambda\gtrsim p^{-1}+q^{-1}]}
\text{\ \ and\ \ } |a'_+(\lambda)|\les \lambda^{-1}(pq)^{-\half}\chi_{[\lambda\gtrsim p^{-1}+q^{-1}]}.
\]
Integrating by parts therefore leads to
\begin{align}
 \eqref{eq:claim4} & \les |t|^{-1} \int_0^\infty \frac{|a_{+}(\lambda)|}{|\phi_+'(\lambda)|^2}\, d\lambda
+ |t|^{-1} \int_0^\infty \frac{|a_{+}'(\lambda)|}{|\phi_+'(\lambda)|}\, d\lambda \nn \\
& \les |t|^{-1} (pq)^{-\half} \int_0^\infty (\lambda+\lambda_0)^{-2}
\chi_{[\lambda\gtrsim p^{-1}+q^{-1}]}\,d\lambda\nn \\
&\qquad + |t|^{-1}(pq)^{-\half} \int_0^\infty \lambda^{-1}(\lambda+\lambda_0)^{-1}
\chi_{[\lambda\gtrsim p^{-1}+q^{-1}]}\,d\lambda \nn \\
& \les |t|^{-1} (pq)^{-\half} \int_0^\infty \lambda^{-2}
\chi_{[\lambda\gtrsim p^{-1}+q^{-1}]}\,d\lambda \les |t|^{-1} \frac{\sqrt{pq}}{p+q} \les |t|^{-1},
\label{eq:sp_2}
\end{align}
and thus~\eqref{eq:claim4} also holds for $\phi_{+}$. We leave the remaining details to the reader.
\end{proof}

We now combine Lemmas~\ref{lem:II_1}, \ref{lem:II_2}, and~\ref{lem:II_3} to
obtain the following lemma. It bounds the contribution of the constant term
in the expansion~\eqref{eq:M_inv_exp}, see also~\eqref{eq:RV_exp}.

\begin{lemma}
\label{lem:QD0Q}
For all test functions $f,g$ and all $t$ one has
\begin{align}
& \left|\int_0^\infty e^{it\lambda^2}\lambda\chi(\lambda) \Bigl\la \big[
R_0^{+}(\lambda^2)v QD_0Q v R_0^{+}(\lambda^2)-R_0^{-}(\lambda^2)v QD_0Q v R_0^{-}(\lambda^2)\big]f,g
\Big\ra \, d\lambda \right| \nn \\
& \les |t|^{-1} \|f\|_1\|g\|_1 \label{eq:adin}
\end{align}
with a constant that only depends on $V$.
\end{lemma}
\begin{proof}
Recall the representation~\eqref{eq:H0} with $H_0^{\pm}(z) = J_0(z)\pm iY_0(z)$.
Hence,
\begin{align}
& R_0^{+}(\lambda^2)(x,x_1)R_0^{+}(\lambda^2)(y_1,y)-R_0^{-}(\lambda^2)(x,x_1)R_0^{-}(\lambda^2)(y_1,y)\nn\\
& = -\frac{i}{8}(Y_0(\lambda|x-x_1|)J_0(\lambda|y-y_1|)+J_0(\lambda|x-x_1|)Y_0(\lambda|y-y_1|)).\label{eq:YoJo}
\end{align}
In addition, we break up the integration region $(0,\infty)$ by means of the partition
 $1=\chi(\lambda)+\tchi(\lambda)$. More precisely, write each resolvent as
\[ R_0^{\pm}(\lambda^2)(x,x_1)=\chi(\lambda|x-x_1|)R_0^{\pm}(\lambda^2)(x,x_1)+
\tchi(\lambda|x-x_1|)R_0^{\pm}(\lambda^2)(x,x_1). \]
This leads to four different terms in~\eqref{eq:adin}. Those terms that contain at least
one $\chi(\lambda|x-x_1|)$ or $\chi(\lambda|y-y_1|)$ we rewrite further using~\eqref{eq:YoJo}.
The other term which involves only $\tchi$ we leave in terms of Hankel functions. Each of these different combinations
is estimated by one of the previous three lemmas.
\end{proof}

Next we turn to the term involving~$S$ in~\eqref{eq:RV_exp}.

\begin{lemma}
\label{lem:S}
Let $S$ and $h_{\pm}(\lambda)$ be as in Lemma~\ref{lem:inv_exp}. Then
for all test functions $f,g$ and all $t$ one has
\begin{align}
& \left|\int_0^\infty e^{it\lambda^2}\lambda\chi(\lambda) \Bigl\la \big[
\frac{1}{h_+(\lambda)}R_0^{+}(\lambda^2)v S v R_0^{+}(\lambda^2)-\frac{1}{h_-(\lambda)}R_0^{-}(\lambda^2)v S v R_0^{-}(\lambda^2)\big]f,g
\Big\ra \, d\lambda \right| \nn \\
& \les |t|^{-1} \|f\|_1\|g\|_1 \label{eq:dwa}
\end{align}
with a constant that only depends on $V$.
\end{lemma}
\begin{proof} Recall that $S$ is of finite rank, and thus Hilbert-Schmidt. In particular,
if $S(x,y)$ denotes the kernel of $S$, then $|S(x,y)|$ is again an $L^2$-bounded operator.
Hence, one shows as before that~\eqref{eq:dwa} reduces to the bound
\begin{align}
&\Big| \int_0^\infty e^{it\lambda^2} \lambda \chi(\lambda) \big[H_0^+(\lambda|x-x_1|)H_0^+(\lambda|y_1-y|)
h_{+}^{-1}(\lambda) \nn \\
&\qquad\qquad\qquad - H_0^-(\lambda|x-x_1|)H_0^-(\lambda|y_1-y|)
h_{-}^{-1}(\lambda) \big] \, d\lambda \Big|  \nn \\
& \les |t|^{-1} (1+\log^{-}|x-x_1|)(1+\log^{-}|y-y_1|).  \label{eq:claim5}
\end{align}
As before, we set $p:=|x-x_1|$ and $q:=|y_1-y|$ for simplicity.
We again need to distinguish whether or not the arguments of the Hankel functions
are $>1$ or $<1$. This will be accomplished by means of the usual partition of
unity $1=\chi+\tchi$. It will also be important to remember that
\[ h_{+}(\lambda)=a\log \lambda+z \text{\ \ and\ \ }h_{-}(\lambda)=a\log \lambda+\bar{z},\]
where $a\ne0$. It is understood that the cut-off $\chi(\lambda)$ in~\eqref{eq:dwa}
is such that $h_{\pm}(\lambda)\ne0$ on the support of $\chi$.
One of the four terms in~\eqref{eq:claim5} which arises as a combination of $\chi$ and $\tchi$
is
\begin{align}
&\Big| \int_0^\infty e^{it\lambda^2} \lambda \chi(\lambda)\chi(\lambda p)\chi(\lambda q)
\frac{J_0(\lambda p)J_0(\lambda q)-Y_0(\lambda p)Y_0(\lambda q)}{(\log\lambda+c_1)^2+c_2^2} \,d\lambda
\Big| \nn \\
&+ \Big| \int_0^\infty e^{it\lambda^2} \lambda \chi(\lambda)\chi(\lambda p)\chi(\lambda q)
\frac{[J_0(\lambda p)Y_0(\lambda q)+Y_0(\lambda p)J_0(\lambda q)](\log\lambda+c_1)}{(\log\lambda+c_1)^2+c_2^2} \,d\lambda
\Big| \nn \\
& \qquad\qquad \les |t|^{-1} (1+\log^{-}p)(1+\log^{-}q) \label{eq:close1}
\end{align}
This is proved by one integration by parts using $\lambda e^{it\lambda^2}=\frac{1}{2it}\frac{d}{d\lambda}e^{it\lambda^2}$.
In view of~\eqref{eq:Y0} the fractions inside of the two integrals take the values $\frac{4}{\pi^2}$
and $\frac{4}{\pi}$, respectively, at $\lambda=0$. Thus, the boundary terms contribute $\les|t|^{-1}$
to the integration by parts. It remains to show that
\begin{align*}
& \int_0^\infty \Big|\frac{d}{d\lambda}\Big[ \chi(\lambda)\chi(\lambda p)\chi(\lambda q)
\frac{J_0(\lambda p)J_0(\lambda q)-Y_0(\lambda p)Y_0(\lambda q)}{(\log\lambda+c_1)^2+c_2^2}\Big]\Big| \,d\lambda
  \\
&+  \int_0^\infty  \Big|\frac{d}{d\lambda}\Big[\chi(\lambda)\chi(\lambda p)\chi(\lambda q)
\frac{[J_0(\lambda p)Y_0(\lambda q)+Y_0(\lambda p)J_0(\lambda q)](\log\lambda+c_1)}{(\log\lambda+c_1)^2+c_2^2} \Big]\Big| \,d\lambda  \\
& \qquad\qquad \les  (1+\log^{-}p)(1+\log^{-}q)
\end{align*}
If the derivative falls on $\chi(\lambda)$, then the resulting term is clearly bounded by
$\les (1+\log^{-}p)(1+\log^{-}q)$. On the other hand, suppose it falls on $\chi(\lambda p)$.
Then that term contributes
\[ \les \frac{1+\log^{-}(q/p)}{1+\log^+p} \chi_{[p\gtrsim 1]} \les 1+\log^{-}q,
\]
and similarly if the derivative falls on~$\chi(\lambda q)$. It therefore remains to check that,
with $\lambda_1=cp^{-1}\wedge cq^{-1}\wedge c$ ($c$ being some small constant)
\begin{align}
& \int_0^{\lambda_1} \Big|\frac{d}{d\lambda}\Big[
\frac{J_0(\lambda p)J_0(\lambda q)-Y_0(\lambda p)Y_0(\lambda q)}{(\log\lambda+c_1)^2+c_2^2}\Big]\Big| \,d\lambda
\label{eq:ARGH1}  \\
&+  \int_0^{\lambda_1} \Big|\frac{d}{d\lambda}\Big[
\frac{[J_0(\lambda p)Y_0(\lambda q)+Y_0(\lambda p)J_0(\lambda q)](\log\lambda+c_1)}{(\log\lambda+c_1)^2+c_2^2} \Big]\Big| \,d\lambda  \label{eq:ARGH2}\\
& \qquad\qquad \les  (1+\log^{-}p)(1+\log^{-}q).\nn
\end{align}
We start with \eqref{eq:ARGH1}. Recall the expansion~\eqref{eq:Y0'} for~$Y_0'$. Also,
let $n(\lambda)>0$ be such that $n(\lambda)^2=(\log\lambda+c_1)^2+c_2^2$.
Then clearly $n(\lambda)\sim|\log\lambda|$ and $n'(\lambda)=\lambda^{-1}+O((\lambda\log\lambda)^{-1})$
as $\lambda\to0$. Hence
\begin{align}
\eqref{eq:ARGH1} &\les \int_0^{\lambda_1} \Big|\frac{\lambda^{-1}+pg(\lambda p)}{n(\lambda)^2}(\log^{-}(\lambda q)+O(1)) + \frac{\lambda^{-1}+qg(\lambda q)}{n(\lambda)^2}(\log^{-}(\lambda p)+O(1)) \nn \\
&\qquad\qquad + 2\frac{(\log^{-}(\lambda p)+O(1))(\log^{-}(\lambda q)+O(1))}{n(\lambda)^3}\,n'(\lambda)
\Big|\,d\lambda + 1. \label{eq:n}
\end{align}
Each of the three terms inside the absolute value contains an expression of the form
$\frac{1}{\lambda \log\lambda}$. Since these are not integrable, one needs to check that they
cancel. Indeed, combining them yields
\begin{equation}
\label{eq:uberf}
 \frac{2\log\lambda}{\lambda n(\lambda)^2}-\frac{2(\log\lambda)^2}{n(\lambda)^3}n'(\lambda)=
O(\lambda^{-1}(\log\lambda)^{-2}),
\end{equation}
which is integrable. Otherwise, we claim that $\eqref{eq:n}\les (1+\log^{-}p)(1+\log^{-}q)$.
To see this, observe first that for all $0<\lambda<\lambda_1$,
\[ \log^{-}(\lambda p) = \log^{-}\lambda+\log^{-}(p)=-\log(\lambda)-\log(p) ,\quad
\log^{-}(\lambda q) =-\log(\lambda)-\log(p).\]
Hence,
\begin{align}
&\Big| \int_0^{\lambda_1} \frac{pg(p\lambda)}{n(\lambda)^2}(\log^{-}(\lambda q) + O(1))\,d\lambda \Big| \nn\\
&\les \int_0^{\lambda_1} \frac{p|g(p\lambda)|}{n(\lambda)^2}(1+\log^{-}(\lambda)) \, d\lambda
(1+|\log q|) \nn \\
&\les \int_0^{p^{-1}} {p|g(p\lambda)|} \, d\lambda
(1+\log^{-} q) \label{eq:pgp} \\
&\les 1+\log^{-} q. \label{eq:I2}
\end{align}
To pass to~\eqref{eq:pgp}, note that if $q\ge1$, then
\[ \sup_{0<\lambda <\lambda_1}\frac{1+\log^{-}(\lambda)}{n(\lambda)^2} (1+\log q)
\les \sup_{0<\lambda < q^{-1}}\frac{1+\log^{-}(\lambda)^2}{n(\lambda)^2}  \les 1, \]
whereas if $0<q<1$, then
\[ \sup_{0<\lambda <\lambda_1}\frac{1+\log^{-}(\lambda)}{n(\lambda)^2} (1+\log^{-} q) \les 1+\log^{-} q. \]
Furthermore,
\begin{align}
&\Big| \int_0^{\lambda_1} \frac{\lambda^{-1}}{n(\lambda)^2}(\log^{-}(\lambda q)-\log\lambda+O(1)) d\lambda
\Big|\nn \\
 &\les \int_0^{cq^{-1}\wedge c} \frac{1}{\lambda (\log\lambda)^2}\,d\lambda\; (1+|\log q|) \nn \\
&\les \frac{1}{\log^{-} (cq^{-1}\wedge c)} (1+|\log q|) \les 1+\log^{-} q. \label{eq:I1}
\end{align}
Finally, we estimate
\begin{align}
& \Big| \int_0^{\lambda_1} \frac{(\log^{-}(\lambda p)+O(1))(\log^{-}(\lambda q)+O(1))-(\log\lambda)^2}{n(\lambda)^3}\,n'(\lambda) \,d\lambda \Big| \nn \\
& \les \int_0^{\lambda_1}\frac{n'(\lambda)}{n(\lambda)^3}\,d\lambda+\int_0^{\lambda_1} \frac{\log^{-}(\lambda q)}{n(\lambda)^3}\,n'(\lambda) \,d\lambda +  \int_0^{\lambda_1} \frac{\log^{-}(\lambda p)}{n(\lambda)^3}\,n'(\lambda) \,d\lambda \label{eq:III1}\\
& + \Big| \int_0^{\lambda_1} \frac{\log^{-}(\lambda p)\log^{-}(\lambda q)-(\log\lambda)^2}{n(\lambda)^3}\,n'(\lambda) \,d\lambda \Big|. \label{eq:III2}
\end{align}
By our previous discussion,
\begin{equation}
\label{eq:III3}
\eqref{eq:III1} \les 1+ \log^{-} q + \log^{-} p.
\end{equation}
On the other hand,
\begin{align}
\eqref{eq:III2} &\les \int_0^{\lambda_1} \frac{|\log p||\log^{-}(\lambda q)| + |\log q|}{\lambda(\log\lambda)^3} \, d\lambda \nn \\
&\les \int_0^{cp^{-1}\wedge c} \frac{|\log p|}{\lambda(\log\lambda)^2} \, d\lambda +
\int_0^{cq^{-1}\wedge cp^{-1}\wedge c } \frac{|\log p| |\log q|}{\lambda (\log \lambda)^3} \, d\lambda \nn\\
& + \int_0^{c q^{-1}\wedge c} \frac{|\log q|}{\lambda(\log\lambda)^3}\, d\lambda \nn \\
&\les (1+\log^{-}p)(1+\log^{-} q).
\label{eq:III4}
\end{align}
Combining \eqref{eq:III4}, \eqref{eq:III3}, \eqref{eq:I2}, \eqref{eq:I1} (and their analogues with
$p$ and $q$ interchanged), as well as~\eqref{eq:uberf}
yields that
\[ \eqref{eq:n} \les (1+\log^{-}p)(1+\log^{-} q),
\]
as claimed.
As far as~\eqref{eq:ARGH2} is concerned, it will suffice to treat the term
involving~$J_0(\lambda p)Y_0(\lambda q)$. This amounts to bounding
\begin{align}
& \int_0^{\lambda_1} \Big|\frac{(\lambda p^2 + O(\lambda^3 p^4))(\log^{-}(\lambda q)+O(1))}{n(\lambda)^2}
(\log\lambda+c_1) \label{eq:J0der}\\
&\qquad +\frac{(1+O(\lambda^2 p^2))(-\lambda^{-1}+qg(q\lambda))}{n(\lambda)^2}(\log\lambda+c_1) \nn\\
&\qquad + \frac{(1+O(\lambda^2 p^2))(\log^{-}(\lambda q)+O(1))}{\lambda n(\lambda)^2} \nn\\
&\qquad -2 \frac{(1+O(\lambda^2 p^2))(\log^{-}(\lambda q)+O(1))}{ n(\lambda)^3}n'(\lambda)(\log\lambda+c_1) \Big|\,d\lambda \label{eq:fsjo}
\end{align}
The first line~\eqref{eq:J0der} contributes $\les 1+\log^{-}q$, as do all the $O$-terms in the
other three lines. The remaining expression inside the absolute values is
\[ -2\frac{\log\lambda}{\lambda n(\lambda)^2}+ 2(\log\lambda)^2\frac{n'(\lambda)}{n(\lambda)^3}=O(\lambda^{-1}(\log\lambda)^{-2})\]
as $\lambda\to0$. This establishes~\eqref{eq:close1}.

Next we turn to the term containing the product $\tchi(\lambda p)\tchi(\lambda q)$.
In analogy with Lemma~\ref{lem:II_3} we work with the Hankel functions rather than $J_0,Y_0$.
Thus we need to show that
\beeq
\label{eq:claim6}
\Bigl| \int_0^\infty e^{i[t\lambda\pm\lambda(p+q)]} \frac{\lambda\chi(\lambda)}{h_{\pm}(\lambda)} \tchi(\lambda p)\tchi(\lambda q) \omega_{\pm}(\lambda q) \omega_{\pm}(\lambda p) \, d\lambda \Bigr| \les |t|^{-1}
\eneq
uniformly in $p,q>0$. Up to the factors $h_{\pm}^{-1}$ this is the same as~\eqref{eq:claim4}.
Combine these factors with the $\lambda$-factor that appears in the integrand. This leads to
functions that satisfy
\[
\Big|\frac{\lambda}{h_{\pm}(\lambda)}\Big| \les \lambda \text{\ \ and\ \ } \Big|\frac{d}{d\lambda}\frac{\lambda}{h_{\pm}(\lambda)}\Big| \les 1
\]
on the support of $\chi$. Hence all the arguments from the proof of Lemma~\ref{lem:II_3} apply to
this case as well, and~\eqref{eq:claim6} holds.

It remains to consider terms that contain $\chi(\lambda p)\tchi(\lambda q)$ or $\chi(\lambda q)\tchi(\lambda p)$. These terms are analogous to those  in Lemma~\ref{lem:II_2}.
We claim that
\begin{align}
&\Big| \int_0^\infty e^{it\lambda^2} \lambda \chi(\lambda)\tchi(\lambda p)\chi(\lambda q)
\frac{J_0(\lambda p)J_0(\lambda q)-Y_0(\lambda p)Y_0(\lambda q)}{(\log\lambda+c_1)^2+c_2^2} \,d\lambda
\Big| \nn \\
&+ \Big| \int_0^\infty e^{it\lambda^2} \lambda \chi(\lambda)\tchi(\lambda p)\chi(\lambda q)
\frac{[J_0(\lambda p)Y_0(\lambda q)+Y_0(\lambda p)J_0(\lambda q)](\log\lambda+c_1)}{(\log\lambda+c_1)^2+c_2^2} \,d\lambda
\Big| \nn \\
& \qquad\qquad \les |t|^{-1} (1+\log^{-}q). \label{eq:mixed}
\end{align}
Write $J_0$, $Y_0$ as
\[
J_0(y)=e^{iy} \rho_+(y)+e^{-iy}\rho_-(y) \text{\ \ and\ \ }Y_0(y)=e^{iy} \sigma_+(y)+e^{-iy}\sigma_-(y)
\]
where $\rho_{\pm}$, $\sigma_{\pm}$ decay like $y^{-\half}$ together with the natural derivative
bounds. Thus \eqref{eq:mixed} is the same as
\begin{align}
&\Big| \int_0^\infty e^{it\psi_{\pm}(\lambda)} \lambda \chi(\lambda)\tchi(\lambda p)\chi(\lambda q)
\frac{\rho_{\pm}(\lambda p)J_0(\lambda q)-\sigma_{\pm}(\lambda p)Y_0(\lambda q)}{(\log\lambda+c_1)^2+c_2^2} \,d\lambda
\Big| \nn \\
&+ \Big| \int_0^\infty e^{it\psi_{\pm}(\lambda)} \lambda \chi(\lambda)\tchi(\lambda p)\chi(\lambda q)
\frac{[\rho_{\pm}(\lambda p)Y_0(\lambda q)+\sigma_{\pm}(\lambda p)J_0(\lambda q)](\log\lambda+c_1)}{(\log\lambda+c_1)^2+c_2^2} \,d\lambda
\Big| \nn \\
& \qquad\qquad \les |t|^{-1} (1+\log^{-}q), \label{eq:mixed2}
\end{align}
where $\psi_{\pm}(\lambda)=\lambda^2\pm \frac{p\lambda}{t}$. The bound~\eqref{eq:mixed2}
can is obtained by means of Lemma~\ref{lem:sp}. In fact, the analysis in Lemma~\ref{lem:II_2}
carries over to this case with minor modifications. To see this, note that
\begin{align}
& \Big|\lambda \chi(\lambda)\tchi(\lambda p)\chi(\lambda q)
\frac{\rho_{\pm}(\lambda p)J_0(\lambda q)-
\sigma_{\pm}(\lambda p)Y_0(\lambda q)}{(\log\lambda+c_1)^2+c_2^2}\Big| \label{eq:mon1}\\
& \les \lambda\chi(\lambda)\tchi(\lambda p)\chi(\lambda q)(\lambda p)^{-\half} (1+\log^{-}q) \nn
\end{align}
and also
\begin{align}
& \Big|\lambda \chi(\lambda)\tchi(\lambda p)\chi(\lambda q)
\frac{[\rho_{\pm}(\lambda p)Y_0(\lambda q)+\sigma_{\pm}(\lambda p)J_0(\lambda q)](\log\lambda+c_1)}{(\log\lambda+c_1)^2+c_2^2}\Big| \label{eq:mon2}\\
& \les \lambda\chi(\lambda)\tchi(\lambda p)\chi(\lambda q)(\lambda p)^{-\half} (1+\log^{-}q). \nn
\end{align}
And similarly for the derivatives. Since these bounds are the same (or even slightly better)
than those satisfied by the functions $a_{\pm}$ in~\eqref{eq:amin} and~\eqref{eq:aplus},
the analysis of Lemma~\ref{lem:II_2} pertaining to these functions carries over to this case as
well, cf.~\eqref{eq:2a}, \eqref{eq:claim3}, and~\eqref{eq:a'again}. This finishes the proof.
\end{proof}

In view of Corollary~\ref{cor:RV_exp}, the only remaining piece in the proof of
Proposition~\ref{prop:low_decay}
is that term in the expansion~\eqref{eq:RV_exp} which involves $E^{\pm}$.

\begin{lemma}
\label{lem:E}
Let $E^{\pm}(\lambda)$ be as in Lemma~\ref{lem:inv_exp}. Then
for all test functions $f,g$ and all $t$ one has
\begin{align}
& \left|\int_0^\infty e^{it\lambda^2}\lambda\chi(\lambda) \Bigl\la \big[
R_0^{+}(\lambda^2)vE^+(\lambda)vR_0^{+}(\lambda^2)-
R_0^{-}(\lambda^2)vE^{-}(\lambda)vR_0^{-}(\lambda^2)\big]f,g
\Big\ra \, d\lambda \right| \nn \\
& \les |t|^{-1} \|f\|_1\|g\|_1 \label{eq:tri}
\end{align}
with a constant that only depends on $V$.
\end{lemma}
\begin{proof}
In analogy with Lemmas~\ref{lem:II_1}, \ref{lem:II_2}, and~\ref{lem:II_3}
we divide the proof into three separate estimates namely,
\begin{align}
& \Big|\int_{\R^8}\int_0^\infty e^{it\lambda^2} \lambda \chi(\lambda) \chi(\lambda|x-x_1|)
H^{\pm}_0(\lambda|x-x_1|) v(x_1) E^{\pm}(\lambda)(x_1,y_1) v(y_1) \nn \\
& H^{\pm}_0(\lambda|y_1-y|)\chi(\lambda|y_1-y|)\,d\lambda\,f(x)g(y)\,dx_1dy_1\,dxdy\Big|
\le C\, |t|^{-1}\|f\|_1\|g\|_1 \label{eq:E_one} \\
& \Big|\int_{\R^8}\int_0^\infty e^{it\lambda^2} \lambda \chi(\lambda) \chi(\lambda|x-x_1|)
H^{\pm}_0(\lambda|x-x_1|) v(x_1)E^{\pm}(\lambda)(x_1,y_1) v(y_1) \nn \\
& H^{\pm}_0(\lambda|y_1-y|)\tchi(\lambda|y_1-y|)\,d\lambda\,f(x)g(y)\,dx_1dy_1\,dxdy\Big|
\le C\, |t|^{-1}\|f\|_1\|g\|_1 \label{eq:E_two} \\
& \Big|\int_{\R^4}\int_0^\infty\int_{\R^4} e^{it\lambda^2} \lambda \chi(\lambda) \tchi(\lambda|x-x_1|) H^{\pm}_0(\lambda|x-x_1|) v(x_1)E^{\pm}(\lambda) (x_1,y_1) v(y_1) \nn \\
& H^{\pm}_0(\lambda|y_1-y|)\tchi(\lambda|y_1-y|)\,d\lambda\,f(x)g(y)\,dx_1dy_1\,dxdy\Big|
\le C\, |t|^{-1}\|f\|_1\|g\|_1 \label{eq:E_three}.
\end{align}
Unlike in the case of $QD_0Q$ we do not exploit any cancellation between $H_0^{+}$
and $H_0^{-}$. This is not only impossible but also unnecessary. In contrast to $QD_0Q$,
the logarithmic singularities of $H_0^{\pm}$ at zero are compensated for by the vanishing
of $E^{\pm}(\lambda)$ at $\lambda=0$, see~\eqref{eq:Eest}.

Let us start with that term where these singularities are not present, i.e., with~\eqref{eq:E_three}.
Set $p=|x-x_1|$, $q=|y-y_1|$, and $\lambda_0=\frac{p+q}{2t}$.
Using the representation~\eqref{eq:ompm} and Lemma~\ref{lem:sp}, we arrive at
\begin{align}
& \Big|\int_0^\infty e^{i[t\lambda^2-\lambda(p+q)]} \lambda \chi(\lambda) \tchi(\lambda p) \omega_{-}(\lambda p) E^{-}(\lambda) (x_1,y_1)
\tchi(\lambda q)\omega_{-}(\lambda q)\,d\lambda \Big| \nn \\
& \les |t|^{-1} \Big\{\int_0^\infty \frac{|a_{-}(\lambda)|}{\delta^2+|\lambda-\lambda_0|^2}\,d\lambda
+ \int_0^\infty \!\!\!\!\chi_{[|\lambda-\lambda_0|>\delta]} \frac{|a_{-}'(\lambda)|}{|\lambda-\lambda_0|}\,d\lambda \Big\}
 \sup_{0<\lambda<\lambda_1}|E^{-}(\lambda) (x_1,y_1)| \nn \\
& + |t|^{-1} \int_0^\infty\!\!\!\! \chi_{[|\lambda-\lambda_0|>\delta]} \frac{\lambda^{-\half}|a_{-}(\lambda)|}{|\lambda-\lambda_0|}\,d\lambda \; \sup_{0<\lambda<\lambda_1}\sqrt{\lambda} |\partial_\lambda
E^{-}(\lambda) (x_1,y_1)|,
\label{eq:a'new}
\end{align}
where we have set
\[ a_{-}(\lambda):= \lambda \chi(\lambda)\tchi(\lambda p) \omega_{-}(\lambda p)
\tchi(\lambda q)\omega_{-}(\lambda q).
\]
Note that the first two integrals involving $a_{-}$ appearing in~\eqref{eq:a'new}
 have already been treated in Lemma~\ref{lem:II_3}. Thus, the expression in braces is $\les 1$.
Moreover, the third integral which involves the new term $\lambda^{-\half}a_{-}(\lambda)$
is actually better than $a_{-}'(\lambda)$, since the latter involves the loss of a full power
of $\lambda$ relative to $a_{-}$ rather than just a half power. Referring to the proof
of Lemma~\ref{lem:II_2} we can therefore again claim that the third integral in~\eqref{eq:a'new}
is $\les 1$. All that remains now is to observe that~\eqref{eq:E_three} follows from the
preceding by means of the error estimates~\eqref{eq:Eest}.

\noindent The case of $E^{+}$ is treated in an analogous fashion, see~\eqref{eq:sp_2},
 and we skip the details.

\noindent Next we consider the other extreme case, i.e., \eqref{eq:E_one} in which $H^{\pm}_0$
is only evaluated on the interval $(0,1]$. Setting
\[  a_{\pm}(\lambda):=  \chi(\lambda)\chi(\lambda p) \omega_{\pm}(\lambda p)
\chi(\lambda q)\omega_{\pm}(\lambda q),
\]
a single integration by parts now yields
\begin{align}
& \Big|\int_0^\infty e^{it\lambda^2} \lambda \chi(\lambda) \chi(\lambda p) \omega_{\pm}(\lambda p)
E^{\pm}(\lambda) (x_1,y_1) \chi(\lambda q)\omega_{\pm}(\lambda q)\,d\lambda \Big| \nn \\
& \les |t|^{-1} \int_0^\infty \sqrt{\lambda}|a_{\pm}'(\lambda)|\,d\lambda
 \sup_{0<\lambda<\lambda_1}\lambda^{-\half}|E^{\pm}(\lambda) (x_1,y_1)| \nn \\
& + |t|^{-1} \int_0^\infty \lambda^{-\half} |a_{\pm}(\lambda)|\,d\lambda\;
\sup_{0<\lambda<\lambda_1}\sqrt{\lambda} |\partial_\lambda E^{\pm}(\lambda) (x_1,y_1)|.
\label{eq:a'new_2}
\end{align}
Now
\begin{align*}
 |a_{\pm}(\lambda)| &\les \chi(\lambda) (1+|\log\lambda|^2)(1+\log^{-} p)(1+\log^{-} q) \\
 |a_{\pm}'(\lambda)| &\les \chi_{[0<\lambda<1]}\lambda^{-1}(1+|\log\lambda|)(1+\log^{-} p)(1+\log^{-} q).
\end{align*}
To obtain \eqref{eq:E_one}, insert these bounds into~\eqref{eq:a'new_2} and invoke~\eqref{eq:Eest}.

\noindent It remains to consider the term of mixed type, i.e., \eqref{eq:E_two}.
Thus set
\[ a_{-}(\lambda):= \lambda \chi(\lambda)\tchi(\lambda p) \omega_{-}(\lambda p)
\chi(\lambda q)\omega_{-}(\lambda q).
\]
Applying Lemma~\ref{lem:sp} with~$\lambda_0=\frac{p}{2t}$ one obtains
\begin{align}
& \Big|\int_0^\infty e^{i[t\lambda^2-\lambda p]} \lambda \chi(\lambda) \tchi(\lambda p) \omega_{-}(\lambda p) E^{-}(\lambda) (x_1,y_1)
\chi(\lambda q)\omega_{-}(\lambda q)\,d\lambda \Big| \nn \\
& \les |t|^{-1} \Big\{\int_0^\infty \frac{\sqrt{\lambda}|a_{-}(\lambda)|}{\delta^2+|\lambda-\lambda_0|^2}\,d\lambda \nn \\
& \qquad\qquad + \int_0^\infty \!\!\!\!\chi_{[|\lambda-\lambda_0|>\delta]} \frac{\sqrt{\lambda}|a_{-}'(\lambda)|}{|\lambda-\lambda_0|}\,d\lambda \Big\}
 \sup_{0<\lambda<\lambda_1}\lambda^{-\half}|E^{-}(\lambda) (x_1,y_1)| \nn \\
& + |t|^{-1} \int_0^\infty\!\!\!\! \chi_{[|\lambda-\lambda_0|>\delta]} \frac{\lambda^{-\half}|a_{-}(\lambda)|}{|\lambda-\lambda_0|}\,d\lambda \; \sup_{0<\lambda<\lambda_1}\sqrt{\lambda} |\partial_\lambda
E^{-}(\lambda) (x_1,y_1)|.
\label{eq:a'new_3}
\end{align}
The basic estimates on $a_{-}(\lambda)$ are
\begin{align*}
|a_{-}(\lambda)| &\les \lambda \chi(\lambda) (\lambda p)^{-\half} (1+\log^{-}\lambda)(1+\log^{-} q) \\
|a_{-}'(\lambda)| &\les \chi_{[0<\lambda<1]} (\lambda p)^{-\half} (1+\log^{-}\lambda)(1+\log^{-} q).
\end{align*}
Hence
\begin{align*}
\sqrt{\lambda} |a_{-}(\lambda)| &\les \lambda \chi(\lambda) (\lambda p)^{-\half}(1+\log^{-} q) \\
\sqrt{\lambda}|a_{-}'(\lambda)|+\lambda^{-\half} |a_{-}(\lambda)| &\les \chi_{[0<\lambda<1]} (\lambda p)^{-\half} (1+\log^{-} q).
\end{align*}
These are precisely the bounds that were used in the proof of
Lemma~\ref{lem:II_2}, and one can therefore repeat the arguments
appearing there, see~\eqref{eq:2a} to~\eqref{eq:a'2}. Finally, the
phase $t\lambda^2+\lambda p$ can be treated as
in~\eqref{eq:a'again}, and we skip the details.
\end{proof}

{\bf Acknowledgement:} The author wishes to thank Monica Visan for
comments on a preliminary version of this paper, as well the
anonymous referee for a very careful reading and many helpful
comments.

\bibliographystyle{amsplain}

\medskip\noindent
\textsc{Division of Astronomy, Mathematics, and Physics, 253-37 Caltech,} \\
\textsc{Pasadena, CA 91125, U.S.A.}\\
{\em email: }\textsf{\bf   schlag@caltech.edu}

\end{document}